# On the Existence of Almost Periodic Solutions with Applications to Global Entrainment

**Iasson Karafyllis[*] and Miroslav Krstic[**]**

[*]Dept. of Mathematics, National Technical University of Athens, Zografou Campus, 15780, Athens, Greece, email: iasonkar@central.ntua.gr; iasonkaraf@gmail.com

[**]Dept. of Mechanical and Aerospace Eng., University of California, San Diego, La Jolla, CA 92093-0411, U.S.A., email: krstic@ucsd.edu

*To Eduardo Sontag on the occasion of his 75th birthday with gratitude for his deep and inspiring ISS and entrainment results*

**Abstract**

This paper provides two results that are useful in the study of the existence and the stability properties of almost periodic solutions for a given dynamical system. The obtained results are generalizations of recent results for periodic systems and are applied to the global entrainment problem in nonlinear time-invariant control systems. It is shown that local exponential stability for the unforced system and input-to-state stability with respect to small inputs can guarantee global entrainment to small almost periodic inputs. In this way, global entrainment is shown in Lotka-Volterra systems with a Volterra-Lyapunov stable interaction matrix. All results can be extended to the uniformly recurrent case.

## 1. Introduction

The problem of existence and attractivity of almost periodic solutions of systems described by ordinary differential equations is an old problem that has attracted the attention of many researchers (see for instance [2, 5, 10] and references therein). It is a difficult mathematical problem which requires special attention because there are counterexamples which show that results in the periodic case are not valid in the almost periodic case.

In this work we are inspired by Eduardo Sontag's results for Input-to-State Stability and entrainment and we continue our recent efforts in [15] to apply conditions that involve matrix measures of the Jacobian matrix in order to explore existence and attractivity of almost periodic solutions. Under almost periodic forcing, a solution that is defined and bounded on the whole real line need not itself be almost periodic; almost periodic systems possessing bounded entire solutions but no almost periodic solutions exist. Here we provide two novel results (Theorem 1 and Theorem 2) that show that a bounded entire solution can inherit the rhythm of the excitation if a fixed amount of contraction can be accumulated within a finite window, uniformly in time and without thinning as time advances (see properties (I) and (II) in Theorem 1 and Theorem 2, respectively). The mechanism compares the solution not with a neighboring trajectory but with its own time-shift by an almost period; recurrent contraction collapses that difference to the residual mismatch of the vector field under the shift, which vanishes along the almost periods. The bounded entire solution is thereby almost periodic, with the same frequency module as the forcing. Theorem 1 is a genuine extension of Theorem 2 in [15]: this is shown by Proposition 1.

Our main results can be applied directly to a well-known problem: the entrainment problem. Entrainment is an important phenomenon arising in many applications ranging from biological systems to the stability of a power grid; see [3, 19, 23] and the discussion in [20]. Indeed, periodic rhythms are present in nature and regulate the function of all living organisms. From a mathematical viewpoint, the problem of proving that entrainment takes place is a very difficult problem. For nonlinear systems, driving the system by an external periodic input does not guarantee that the solution of the system converges to a periodic solution, as more complicated phenomena can happen (harmonic generation or suppression, chaos or quasiperiodic solutions). This may happen even if the system is well-behaved in the input-free case; see [23] and [4].

However, in nature most rhythms are not exactly periodic but rather almost periodic because the rhythms come from the superposition of independent periodic phenomena. Therefore, it makes sense to consider the problem of entrainment to almost periodic inputs. There are only few works in the literature that study entrainment to almost periodic inputs. The problem of entrainment to almost periodic inputs was recently studied for the class of Lur'e systems in [7, 8, 9]. Since physical and biological systems usually can tolerate excitation only of bounded magnitude, we frame entrainment in terms of inputs of small magnitude (see Definition 3 below), which we regard as the appropriate notion for such systems.

The present work provides a verifiable entrainment criterion (Corollary 1) requiring only two ingredients: uniform ultimate boundedness of trajectories under small inputs, and a single matrix inequality imposed on a bounded region rather than on the whole state space. The latter is a Demidovich/convergent-systems condition (see [18, 22]); the point of our formulation is that contraction is demanded only where the trajectories ultimately reside, with boundedness supplied separately. This localization is what makes the criterion usable, since the global form of the condition is frequently too strong. Under these two conditions every sufficiently small almost periodic (or periodic) excitation produces a unique response of the same frequency module to which all trajectories converge. A qualitative consequence (Corollary 2) is worth isolating, because it requires no structure beyond properties that are routinely verified. If the system is input-to-state stable with respect to small inputs (see [21]), and the unforced origin is locally exponentially stable, then the system entrains to all sufficiently small almost periodic and periodic excitations. In short, input-to-state stability together with an exponentially stable linearization implies entrainment to small almost periodic disturbances.

Finally, we apply the obtained results to Lotka–Volterra population models with external forcing (Theorem 3). Lotka-Volterra models are popular because of their simplicity and their efficiency to describe qualitatively various phenomena in ecology. Many researchers have studied dynamical properties of Lotka-Volterra models, including boundedness of solutions, persistence, permanence, existence and stability of equilibria, limit cycles and heteroclinic orbits; see for instance [12]. For a system possessing a coexistence equilibrium whose interaction matrix is Volterra–Lyapunov stable we show that small almost periodic or periodic forcing (seasonal hunting, fluctuating resource supply, and superpositions of further environmental cycles) drives every trajectory to a unique almost periodic population rhythm carrying exactly the frequencies of the forcing. Lotka-Volterra systems with a Volterra–Lyapunov stable interaction matrix allow the construction of Lyapunov functions (see the discussions in [17] and [14]).

The structure of paper is as follows. First, we give the basic notions and the notation that are used in the paper. Section 2 provides the main results (Theorem 1 and Theorem 2) of the present work.

Section 3 presents two examples that illustrate the applicability of the main results. The main results (Corollary 1 and Corollary 2) for the entrainment problem are given in Section 4. Section 5 is devoted to the study of the entrainment for the class of Lotka-Volterra systems. Section 6 provides the proofs of all results. Section 7 of the paper deals with the extension of the obtained results to the uniformly recurrent case. The concluding remarks of this work are given in Section 8. Finally, the derivation of an estimate that is used in one of the examples of Section 3 is presented in detail in the Appendix.

**Notation and Basic Notions**

In this paper we adopt the following notation.

* $\mathbb{R}_+ := [0,+\infty)$.
* Let $D \subseteq \mathbb{R}^n$ be a non-empty open set and let $S \subseteq \mathbb{R}^n$ be a set that satisfies $D \subseteq S \subseteq \bar{D}$. By $C^0(S\,;\Omega)$, we denote the class of continuous functions on $S$, which take values in $\Omega \subseteq \mathbb{R}^m$. By $C^k(S\,;\Omega)$, where $k \geq 1$ is an integer, we denote the class of functions on $S \subseteq \mathbb{R}^n$, which take values in $\Omega \subseteq \mathbb{R}^m$ and have continuous derivatives of order $k$. In other words, the functions of class $C^k(S\,;\Omega)$ are the functions which have continuous derivatives of order $k$ in $D$ that can be continued continuously to all points in $S$. When $\Omega = \mathbb{R}$ then we write $C^0(S)$ or $C^k(S)$.
* For $x \in \mathbb{R}$ we define the integer part of $x$ denoted by $\lfloor x \rfloor$ as the largest integer that is less than or equal to $x$, i.e., $\lfloor x \rfloor = \max\{\, q \in \mathbb{Z} : q \leq x \,\}$.
* For a bounded function $x : \mathbb{R} \to \mathbb{R}^n$ we define $\|x_\infty\| = \sup_{t\in\mathbb{R}}\left(|x(t)|\right)$, where $|\cdot|$ denotes the standard Euclidean norm of $\mathbb{R}^n$.
* The shift operator $\delta_\tau$, where $\tau \in \mathbb{R}$ is a constant, is the operator that maps a function $x : \mathbb{R} \to \mathbb{R}^n$ to the "shifted" function $\delta_\tau x : \mathbb{R} \to \mathbb{R}^n$ defined by $(\delta_\tau x)(t) = x(t+\tau)$ for all $t \in \mathbb{R}$.
* By $K$ we denote the class of increasing, continuous functions $a : \mathbb{R}_+ \to \mathbb{R}_+$ with $a(0) = 0$. By $KL$ we denote the set of all continuous functions $\sigma : \mathbb{R}_+ \times \mathbb{R}_+ \to \mathbb{R}_+$ with the properties: (i) for each $t \geq 0$ the mapping $\sigma(\cdot,t)$ is of class $K$; (ii) for each $s \geq 0$, the mapping $\sigma(s,\cdot)$ is non-increasing with $\lim_{t\to+\infty}\left(\sigma(s,t)\right) = 0$.

In the literature there are various versions of almost periodic functions. Here we adopt the notion of almost periodic functions that is used in [10] and [5].

**Definition 1:** *A continuous function* $\varphi : \mathbb{R} \to \mathbb{R}^n$ *is called almost periodic if for every* $\varepsilon > 0$ *there exists* $T(\varepsilon) > 0$ *such that for every* $a \in \mathbb{R}$ *there exists* $\tau \in [a, a+T(\varepsilon)]$ *with* $\|\delta_\tau \varphi - \varphi\|_\infty \leq \varepsilon$.

**Definition 2:** *A continuous map* $f : \mathbb{R} \times \mathbb{R}^n \to \mathbb{R}^m$ *is called almost periodic in* $t$ *uniformly in* $x \in \Lambda$, *where* $\Lambda \subseteq \mathbb{R}^n$ *is a non-empty set, if for every* $\varepsilon > 0$ *there exists* $T(\varepsilon) > 0$ *such that for every* $a \in \mathbb{R}$ *there exists* $\tau \in [a, a+T(\varepsilon)]$ *with* $\sup_{t\in\mathbb{R}, x\in\Lambda}\left(|f(t+\tau,x) - f(t,x)|\right) \leq \varepsilon$.

The following facts about almost periodic functions can be found in [10].

**Fact I:** *Every map* $f:\mathbb{R}\times\mathbb{R}^n\to\mathbb{R}^m$ *which is almost periodic in* $t$ *uniformly in* $x\in\Lambda$, *where* $\Lambda\subseteq\mathbb{R}^n$ *is a non-empty set, has a set* $\Omega\subset\mathbb{R}$ *of basic frequencies, i.e.,* $\Omega\subset\mathbb{R}$ *is a countable set which is linearly independent over the rationals (i.e., if* $\sum_p b_p\omega_p=0$ *with* $\{b_p\}\subseteq\mathbb{Q}$ *and* $\{\omega_p\}\subseteq\Omega$ *then* $b_p=0$*) and for every sequence* $\{\tau_k:k=1,2,...\}\subseteq\mathbb{R}$ *with* $\lim_{k\to+\infty}\left(\exp\left(i\omega_p\tau_k\right)\right)=1$ *for each* $\omega_p\in\Omega$ *it holds that* $\lim_{k\to+\infty}\left(\sup_{t\in\mathbb{R},x\in\Lambda}\left(\left|f(t+\tau_k,x)-f(t,x)\right|\right)\right)=0$.

**Fact II:** *If* $f:\mathbb{R}\times\mathbb{R}^n\to\mathbb{R}^m$ *is a continuous map for which there exists a countable set* $\Omega\subset\mathbb{R}$ *which is linearly independent over the rationals (i.e., if* $\sum_p b_p\omega_p=0$ *with* $\{b_p\}\subseteq\mathbb{Q}$ *and* $\{\omega_p\}\subseteq\Omega$ *then* $b_p=0$*) and if* $\Lambda\subseteq\mathbb{R}^n$ *is a non-empty set with the property that for every sequence* $\{\tau_k:k=1,2,...\}\subseteq\mathbb{R}$ *with* $\lim_{k\to+\infty}\left(\exp\left(i\omega_p\tau_k\right)\right)=1$ *for each* $\omega_p\in\Omega$ *it holds that* $\lim_{k\to+\infty}\left(\sup_{t\in\mathbb{R},x\in\Lambda}\left(\left|f(t+\tau_k,x)-f(t,x)\right|\right)\right)=0$ *then* $f$ *is almost periodic in* $t$ *uniformly in* $x\in\Lambda$ *and* $\Omega\subset\mathbb{R}$ *is the set of basic frequencies of* $f$.

We notice that if $f:\mathbb{R}\times\mathbb{R}^n\to\mathbb{R}^m$ is almost periodic in $t$ uniformly in $x\in\Lambda$ and its set $\Omega\subset\mathbb{R}$ of basic frequencies is a singleton then $f$ restricted on $\mathbb{R}\times\Lambda$ is a $T-$periodic map for some $T>0$.

Let $|\cdot|$ denote a norm of $\mathbb{R}^n$. The induced matrix norm $\|\cdot\|$ is defined for a matrix $G\in\mathbb{R}^{n\times n}$ by the formula $\|G\|:=\max\left\{|Gx|:|x|\le 1\right\}=\max\left\{|Gx|:|x|=1\right\}$. The associated matrix measure or logarithmic norm $\mu(G)$ of $G\in\mathbb{R}^{n\times n}$ is defined by

$$\mu(G):=\lim_{h\to 0^+}\left(\frac{\|I_n+hG\|-1}{h}\right) \tag{1.1}$$

where $I_n$ is the unit $n\times n$ matrix. The book [24] presents many properties of the matrix measure. For the weighted 2-norm $|x|=|Lx|_2$, where $|\cdot|_2$ denotes the standard Euclidean norm in $\mathbb{R}^n$ and $L\in\mathbb{R}^{n\times n}$ is an invertible matrix, we get $\max\left\{y^T\left(PG+G^TP\right)y:y^TPy=1\right\}=2\mu(G)$, where $^T$ denotes transposition and $P=L^TL$. For the $\infty$-norm $|x|=|x|_\infty=\max_{i=1,...,n}\left(|x_i|\right)$ we get $\mu(G)=\max_{i=1,...,n}\left(G_{i,i}+\sum_{j\ne i}\left|G_{i,j}\right|\right)$ and for the 1-norm $|x|=|x|_1=\sum_{i=1}^n|x_i|$ we get $\mu(G)=\max_{j=1,...,n}\left(G_{j,j}+\sum_{i\ne j}\left|G_{i,j}\right|\right)$.

## 2. Main Results

Let $f \in C^1\left(\mathbb{R}\times\mathbb{R}^n;\mathbb{R}^n\right)$ be given. Consider the time-varying system given by the differential equation

$$\dot{x}(t) = f\left(t, x(t)\right) \tag{2.1}$$

where $x(t) \in \mathbb{R}^n$ is the state variable.

Let $t_0 \in \mathbb{R}$, $x_0 \in \mathbb{R}^n$ be given. By $\varphi(t,t_0,x_0)$ we denote the unique solution of (2.1) at time $t \in \mathbb{R}$ with initial condition $x(t_0)=x_0$. For each $t_0 \in \mathbb{R}$, $x_0 \in \mathbb{R}^n$ there exist $t_{\max} \in \left(t_0,+\infty\right]$ and $t_{\min} \in \left[-\infty,t_0\right)$ for which the solution $\varphi(t,t_0,x_0)$ is defined for $t \in \left(t_{\min},t_{\max}\right)$. Moreover, if $t_{\max} < +\infty$ then $\limsup_{t\to t_{\max}^-}\left(\left|\varphi(t,t_0,x_0)\right|\right)=+\infty$, where $\left|\cdot\right|$ denotes a norm of $\mathbb{R}^n$. Similarly, if $t_{\min} > -\infty$ then $\limsup_{t\to t_{\min}^+}\left(\left|\varphi(t,t_0,x_0)\right|\right)=+\infty$. The map $\varphi$ satisfies $\varphi(t_0,t_0,x_0)=x_0$ (identity property) and $\varphi(t,t_0,x_0)=\varphi\left(t,\tau,\varphi(\tau,t_0,x_0)\right)$ (semigroup property) for all $t_0 \in \mathbb{R}$, $x_0 \in \mathbb{R}^n$ and $t,\tau \in \left(t_{\min},t_{\max}\right)$.

We use the term "entire solution of (2.1)" when $\varphi(t,t_0,x_0)$ is defined for all $t \in \mathbb{R}$.

We next make the following assumption for system (2.1).

**(A)** *For every non-empty compact set* $S \subset \mathbb{R}^n$ *it holds that* $\sup\left\{\left|f(t,x)\right| : t\in\mathbb{R}, x\in S\right\} < +\infty$ *and there exists* $L>0$ *such that* $\left|f(t,x)-f(t,y)\right| \le L\left|x-y\right|$ *for all* $t\in\mathbb{R}$, $x,y\in S$.

Recently, in [15] we gave the following result that relates matrix measures and periodic systems.

**Theorem 2 in [15]-Existence of an attracting periodic solution:** *Let* $f \in C^1\left(\mathbb{R}\times\mathbb{R}^n;\mathbb{R}^n\right)$ *be a* $T-$*periodic map for some* $T>0$, *i.e.* $f\left(t+T,x\right)=f\left(t,x\right)$ *for all* $t\in\mathbb{R}$ *and* $x\in\mathbb{R}^n$. *Consider system (2.1) and suppose that there exists a non-empty set* $S\subseteq\mathbb{R}^n$ *and a convex, compact set* $C\subseteq\mathbb{R}^n$ *such that for all* $t_0\ge 0$, $x_0\in S$, *there exists* $\tau\ge t_0$ *with* $\varphi\left(t,t_0,x_0\right)\in C$ *for all* $t\ge\tau$. *Moreover, suppose that there exists a* $T-$*periodic function* $p\in C^0\left(\mathbb{R}\right)$ *with* $\int_0^T p\left(s\right)ds>0$ *such that* $\mu\left(\frac{\partial f}{\partial z}(t,z)\right)\le -p(t)$ *for all* $t\in\mathbb{R}$, $z\in C$, *where* $\mu\left(\cdot\right)$ *denotes the matrix measure associated with a norm* $\left|\cdot\right|$ *of* $\mathbb{R}^n$. *Then there exists a unique* $T-$*periodic solution* $x^*\in C^1\left(\mathbb{R};C\right)$ *of (2.1) such that* $\lim_{t\to+\infty}\left(\left|\varphi\left(t,t_0,x_0\right)-x^*(t)\right|\right)=0$ *for all* $t_0\ge 0$, $x_0\in S$.

A generalization of Theorem 2 in [15] to the almost periodic case is not straightforward. In the almost periodic case there are systems with entire bounded solutions and no almost periodic solutions (see [5] and references therein). Additional properties must be assumed about $f \in C^1\left(\mathbb{R}\times\mathbb{R}^n\right)$ (e.g., $n=1$ and monotonicity of $f$ with respect to $x$; see Theorem 12.8 on page 229 in [5]).

Here we give conditions that involve matrix measures. Our main result is stated next.

**Theorem 1-Existence of an attracting almost periodic solution:** *Consider system (2.1) and assume that* $f \in C^1\left(\mathbb{R}\times\mathbb{R}^n;\mathbb{R}^n\right)$ *satisfies assumption (A). Let* $C \subseteq \mathbb{R}^n$ *be a convex, compact set for which* $f \in C^1\left(\mathbb{R}\times\mathbb{R}^n;\mathbb{R}^n\right)$ *is almost periodic in* $t$ *uniformly in* $x \in C$. *Suppose that there exists* $t_0 \in \mathbb{R}$, $x_0 \in \mathbb{R}^n$ *and* $\tau \in \mathbb{R}$ *with* $\varphi\left(t,t_0,x_0\right) \in C$ *for all* $t \geq \tau$. *Let* $p \in C^0\left(\mathbb{R}\right)$ *be a function that satisfies* $\int_0^{+\infty} p\left(s\right)ds = +\infty$ *and* $\mu\left(\frac{\partial f}{\partial z}\left(t,z\right)\right) \leq -p\left(t\right)$ *for all* $t \in \mathbb{R}$, $z \in C$, *where* $\mu\left(\,\cdot\,\right)$ *denotes the matrix measure associated with a norm* $\left|\cdot\right|$ *of* $\mathbb{R}^n$. *Suppose that there exist constants* $K,c>0$ *such that for every* $t \in \mathbb{R}$ *the following property holds:*

**(I)** *there exists* $t_0 \in \left(-\infty,t\right)$ *with* $\int_{t_0}^{t} p(s)ds \geq c$ *and* $\int_{t_0}^{t} \exp\left(-\int_{l}^{t} p(s)ds\right)dl \leq K$,

*Then there exists an almost periodic solution* $x^* \in C^1\left(\mathbb{R};C\right)$ *of (2.1) which has the same basic frequencies with* $f$. *Moreover,* $\lim_{t\to+\infty}\left(\left|\varphi\left(t,t_0,x_0\right)-x^*(t)\right|\right)=0$ *for all* $t_0 \in \mathbb{R}$, $x_0 \in S\left(t_0\right)$, *where* $S\left(t_0\right) \subseteq \mathbb{R}^n$ *is the set of all* $x_0 \in \mathbb{R}^n$ *for which there exists* $\tau \in \mathbb{R}$ *with* $\varphi\left(t,t_0,x_0\right) \in C$ *for all* $t \geq \tau$.

Theorem 3 is a genuine extension of Theorem 2 in [15]. This is shown by the following proposition.

**Proposition 1:** *Let* $p \in C^0\left(\mathbb{R}\right)$ *be a* $T-$*periodic function with* $\int_0^{T} p\left(s\right)ds > 0$. *Then* $\int_0^{+\infty} p\left(s\right)ds = +\infty$ *and there exist constants* $K,c>0$ *such that for every* $t \in \mathbb{R}$ *property (I) of Theorem 1 holds.*

The condition $\mu\left(\frac{\partial f}{\partial z}\left(t,z\right)\right) \leq -p\left(t\right)$ for all $t \in \mathbb{R}$, $z \in C$ resembles the usual contraction condition $\mu\left(\frac{\partial f}{\partial z}\left(t,z\right)\right) \leq -s$, where $s>0$ is a constant (see [19, 20]) and allows the interpretation of $p(t)$ as the instantaneous contraction rate in $C \subseteq \mathbb{R}^n$. In this way, the asymptotic condition $\int_0^{+\infty} p\left(s\right)ds = +\infty$ can be interpreted as an asymptotic contraction property in the convex, compact set $C \subseteq \mathbb{R}^n$. Moreover, using this interpretation, it can be said that property (I) requires, uniformly in $t$, that a fixed contraction budget $c>0$ accumulates over some finite backward window whose discounted length stays bounded by $K>0$; this is a sort of a recurrent, non-thinning contraction.

The proof of Theorem 1 is based on the following sequence of steps:

> First, we obtain an entire bounded solution → then we show that the bounded entire solution is almost periodic → then we show that the entire bounded solution is attracting → and finally, Theorem 1 assembles everything.

The first step in the roadmap above is similar to see Theorem 6.2 on page 99 in [5]. For reasons of completeness, we next give a result similar to Theorem 6.2 on page 99 in [5].

**Lemma 1-Conditions for the existence of a bounded entire solution:** *Let $C \subseteq \mathbb{R}^n$ be a compact set for which there exists $t_0 \in \mathbb{R}$, $x_0 \in \mathbb{R}^n$ and $\tau \in \mathbb{R}$ with $\varphi(t,t_0,x_0) \in C$ for all $t \ge \tau$. Suppose that $f \in C^1\left(\mathbb{R}\times\mathbb{R}^n;\mathbb{R}^n\right)$ is almost periodic in $t$ uniformly in $x \in C$ and satisfies assumption (A). Then there exists an entire solution $x^* \in C^1(\mathbb{R};C)$ of (2.1).*

The proof of Lemma 1 is provided in Section 6 for reasons of completeness and the proof is slightly different from the proof in [5]. Moreover, the proof of Lemma 1 that is given here allows the extension of Lemma 1 to the uniformly recurrent case: see Section 7 below.

The second step in our roadmap is a novel result that guarantees almost periodicity for entire bounded solutions.

**Theorem 2-Conditions for a bounded entire solution to be an almost periodic solution:** *Consider system (2.1) and let $C \subseteq \mathbb{R}^n$ be a convex, compact set for which $f \in C^1\left(\mathbb{R}\times\mathbb{R}^n;\mathbb{R}^n\right)$ is almost periodic in $t$ uniformly in $x \in C$. Let $p,\beta \in C^0(\mathbb{R})$ be functions that satisfy*

$$\mu\left(\frac{\partial f}{\partial z}(t,z)\right) \le -p(t) \ , \ \mu\left(-\frac{\partial f}{\partial z}(t,z)\right) \le -\beta(t) \tag{2.2}$$

*for all $t \in \mathbb{R}$, $z \in C$, where $\mu(\cdot)$ denotes the matrix measure associated with a norm $|\cdot|$ of $\mathbb{R}^n$. Suppose that there exist constants $K,c>0$ such that for every $t \in \mathbb{R}$ either property (I) of Theorem 1 holds or the following property holds:*

**(II)** *there exists $t_0 \in (t,+\infty)$ with $\int_t^{t_0} \beta(s)ds \ge c$ and $\int_t^{t_0} \exp\left(-\int_t^{l} \beta(s)ds\right)dl \le K$.*

*Furthermore, suppose that assumption (A) holds. Let $x^* : \mathbb{R} \to C$ be an entire solution of (2.1). Then $x^*$ is almost periodic and $x^*$ has the same basic frequencies with $f$.*

The proof of Theorem 2 is given in Section 6 below. Property (II) in Theorem 2 is nothing but property (I) of Theorem 1 for the time-reversed system (2.1). The reason that a different instantaneous contraction rate is used for the time-reversed system (2.1) ( $p(t)$ is the instantaneous contraction rate in forward time and $\beta(t)$ is the instantaneous contraction rate in backward time) is the fact that for a matrix $G \in \mathbb{R}^{n\times n}$ we have $\mu(-G) \neq -\mu(G)$ (in general). Finally, it should be noticed that the function $\beta \in C^0(\mathbb{R})$ involved in (2.2) is not used in Theorem 1: this happens because Theorem 1 deals with the attractivity properties in forward time and does not deal with backward time. On the other hand, Theorem 2 does not deal with attractivity properties of the almost periodic solution and that is why we exploit both time directions.

The third step in the aforementioned sequence of steps is to give conditions that guarantee the attractivity of an entire bounded solution. This is achieved by the following result.

**Lemma 2- Conditions for a bounded entire solution to attract other solutions:** *Let $C \subseteq \mathbb{R}^n$ be a convex, compact set and let $p \in C^0(\mathbb{R})$ be a function that satisfies $\mu\left(\frac{\partial f}{\partial z}(t,z)\right) \leq -p(t)$ for all $t \in \mathbb{R}$, $z \in C$, where $\mu(\cdot)$ denotes the matrix measure associated with a norm $|\cdot|$ of $\mathbb{R}^n$. Suppose that $\int_0^{+\infty} p(s)ds = +\infty$. Let $x^* : \mathbb{R} \to C$ be an entire solution of (2.1). Define for every $t_0 \in \mathbb{R}$ the non-empty set $S(t_0) \subseteq \mathbb{R}^n$ to be the set for which for all $x_0 \in S(t_0)$ there exists $\tau \in \mathbb{R}$ with $\varphi(t,t_0,x_0) \in C$ for all $t \geq \tau$. Then $\lim_{t\to+\infty}\left(\left|\varphi(t,t_0,x_0) - x^*(t)\right|\right) = 0$ for all $t_0 \in \mathbb{R}$, $x_0 \in S(t_0)$.*

The proof of Lemma 2 is given in Section 6 below. Again, it should be noticed that the function $\beta \in C^0(\mathbb{R})$ involved in (2.2) is not used in Lemma 2: this happens because Lemma 2 deals with the attractivity properties in forward time and does not deal with backward time.

## 3. Examples

In this section we provide some examples that illustrate the significance of Theorem 1.

**Example 1:** Consider the scalar system

$$\dot{x}(t) = \left(\kappa - \sin^2(t)\right)x(t) - u(t)x^3(t) + Rh(t) \tag{3.1}$$

where $\kappa \in (0,1/2)$ and $R \geq 0$ are constants, while $u,h \in C^1(\mathbb{R})$ are almost periodic functions that satisfy the following inequalities for all $t \in \mathbb{R}$:

$$u(t) \geq 0 \text{ and } |h(t)| \leq 1 \tag{3.2}$$

It is clear that system (3.1) is a system of the form (2.1) with $f(t,x) = \left(\kappa - \sin^2(t)\right)x - u(t)x^3 + Rh(t)$ for all $t \in \mathbb{R}, x \in \mathbb{R}$. Clearly, $f \in C^1(\mathbb{R}\times\mathbb{R})$ satisfies assumption (A) and for every non-empty compact set $\Lambda \subseteq \mathbb{R}$, $f:\mathbb{R}\times\mathbb{R}\to\mathbb{R}$ is almost periodic in $t$ uniformly in $x\in\Lambda$.

For the scalar system (3.1) we cannot apply Theorem 12.8 on page 229 in [5] since $f$ is not monotone with respect to $x$. To see this, notice that $\frac{\partial f}{\partial x}(t,x)=\kappa-\sin^2(t)-3u(t)x^2$ for all $t\in\mathbb{R}, x\in\mathbb{R}$ with $\frac{\partial f}{\partial x}\left((2j+1)\frac{\pi}{2},x\right)=\kappa-1-3u\left((2j+1)\frac{\pi}{2}\right)x^2\le\kappa-1<0$ for all $j\in\mathbb{Z}, x\in\mathbb{R}$ (a consequence of (3.2) and the fact that $\kappa<1/2$) and $\frac{\partial f}{\partial x}(j\pi,0)=\kappa>0$ for all $j\in\mathbb{Z}$. We cannot also see how the results in [10] can be applied even if we have additional assumptions on $u,h\in C^1(\mathbb{R})$ (e.g. $|h(t)|\le\varepsilon$ and $|u(t)|\le\varepsilon$ for a sufficiently small $\varepsilon>0$) and even then the results in [10] would give local results.

We next show how easily Theorem 1 can be applied to system (3.1) to obtain global results. For this example, due to (3.2) it holds that

$$\mu\left(\frac{\partial f}{\partial x}(t,x)\right)=\kappa-\sin^2(t)-3u(t)x^2\le\kappa-\sin^2(t) \tag{3.3}$$

for all $t\in\mathbb{R}, x\in\mathbb{R}$. Let arbitrary $t_0\in\mathbb{R}$, $x_0\in\mathbb{R}$ be given. Then there exist $t_{\max}\in(t_0,+\infty]$ and $t_{\min}\in[-\infty,t_0)$ for which the solution $\varphi(t,t_0,x_0)$ of (3.1) is defined for $t\in(t_{\min},t_{\max})$. The following implicit solution formula is valid for all $t\in[t_0,t_{\max})$:

$$\begin{aligned}\varphi(t,t_0,x_0)&=\exp\left(-\left(\frac{1}{2}-\kappa\right)(t-t_0)+\frac{1}{4}\left(\sin(2t)-\sin(2t_0)\right)-\int_{t_0}^{t}u(l)\varphi^2(l,t_0,x_0)\,ds\right)x_0\\&+R\int_{t_0}^{t}h(s)\exp\left(-\left(\frac{1}{2}-\kappa\right)(t-s)+\frac{1}{4}\left(\sin(2t)-\sin(2s)\right)-\int_{s}^{t}u(l)\varphi^2(l,t_0,x_0)\,ds\right)ds\end{aligned} \tag{3.4}$$

Using (3.2) and (3.4) we conclude and the following solution estimate holds for all $t\in[t_0,t_{\max})$:

$$|\varphi(t,t_0,x_0)|\le\exp\left(-\left(\frac{1}{2}-\kappa\right)(t-t_0)+\frac{1}{2}\right)|x_0|+\frac{2R}{1-2\kappa}\exp\left(\frac{1}{2}\right) \tag{3.5}$$

Estimate (3.5) shows that $t_{\max}=+\infty$. Moreover, estimate (3.5) shows that $\limsup_{t\to+\infty}\left(|\varphi(t,t_0,x_0)|\right)<b$, where $b>\frac{2R}{1-2\kappa}\exp\left(\frac{1}{2}\right)$ is an arbitrary constant. This fact in conjunction with inequality (3.3), Proposition 1 and the fact that $\int_0^{2\pi}\left(\sin^2(s)-\kappa\right)ds=\pi-2\kappa\pi>0$ imply that all assumptions of Theorem 1 are valid with $S(t_0)\equiv\mathbb{R}$, $p(t):=\sin^2(t)-\kappa$ and $C=[-b,b]$.

It follows from Theorem 1 that there exists an almost periodic solution $x^* \in C^1(\mathbb{R};C)$ of (3.1) which has the same basic frequencies with $f(t,x) = (\kappa - \sin^2(t))x - u(t)x^3 + Rh(t)$. Moreover, $\lim_{t\to+\infty}\left(\left|\varphi(t,t_0,x_0) - x^*(t)\right|\right) = 0$ for all $t_0 \in \mathbb{R}$, $x_0 \in \mathbb{R}$. $\lhd$

The following example shows that Theorem 1 can be applied to time-invariant systems in order to show existence/attractivity of almost periodic solutions.

**Example 2:** Consider the following nonlinear system

$$
\begin{aligned}
\dot{z}_1 &= z_1 - z_2 - \left(z_1^2 + z_2^2\right) z_1 + a z_1 z_3 \\
\dot{z}_2 &= z_1 + z_2 - \left(z_1^2 + z_2^2\right) z_2 + a z_2 z_3 \\
\dot{z}_3 &= z_3 - \omega z_4 - \left(z_3^2 + z_4^2\right) z_3 + b z_1 z_3 \\
\dot{z}_4 &= \omega z_3 + z_4 - \left(z_3^2 + z_4^2\right) z_4 + b z_1 z_4 \\
z &= \left(z_1, z_2, z_3, z_4\right) \in \mathbb{R}^4
\end{aligned}
\tag{3.6}
$$

where $\omega > 0$, $a,b \in \mathbb{R}$ are constants. We show next by using Theorem 1 that under the assumption that the constants $a,b \in \mathbb{R}$ are sufficiently small and so that

$$
\begin{aligned}
&4a^2 + 4b^2 + |a| + |b| < 1 \\
&\max\left(|a|,|b|\right)\exp\left(\kappa\left(2\sqrt{\frac{a^2+b^2}{1-|a|-|b|}}\right)\right) < 2
\end{aligned}
\tag{3.7}
$$

where $\kappa(s) = s + \sqrt{2s}$ for $s \geq 0$, then every solution of system (3.6) tends asymptotically to an almost periodic solution of (3.6).

System (3.6) in polar coordinates

$$
\begin{aligned}
z_1 &= r\cos(\theta) \quad , \quad z_2 = r\sin(\theta) \\
z_3 &= \rho\cos(\varphi) \quad , \quad z_4 = \rho\sin(\varphi)
\end{aligned}
\tag{3.8}
$$

is described by the differential equations

$$
\begin{aligned}
&\dot{r} = r\left(1 - r^2 + a\rho\cos(\varphi)\right) \quad , \quad \dot{\rho} = \rho\left(1 - \rho^2 + br\cos(\theta)\right) \\
&\dot{\theta} = 1 \quad , \quad \dot{\varphi} = \omega
\end{aligned}
\tag{3.9}
$$

Since $\theta(t) = \theta_0 + t$, $\varphi(t) = \varphi_0 + \omega t$ for some $\theta_0, \varphi_0 \in \mathbb{R}$ we can look for solutions of (3.9) which satisfy the time-varying differential equations:

$$
\begin{aligned}
\dot{r}(t) &= r(t)\left(1 - r^2(t) + a\rho(t)\cos\left(\varphi_0 + \omega t\right)\right) \\
\dot{\rho}(t) &= \rho(t)\left(1 - \rho^2(t) + br(t)\cos\left(\theta_0 + t\right)\right)
\end{aligned}
\tag{3.10}
$$

The study of system (3.10) is easy when $r(t) \equiv 0$ or $\rho(t) \equiv 0$. When $\rho(t) \equiv 0$ (i.e., when $z_3(0) = z_4(0) = 0$ in the original coordinates of (3.6)) then the solution of system (3.10) either satisfies $r(t)\equiv 0$ or $\lim_{t\to+\infty}\left(r(t)\right)=1$. In this case the solution of (3.6) either tends to an equilibrium point or tends to a periodic solution of period $2\pi$. When $r(t)\equiv 0$ (i.e., when $z_1(0)=z_2(0)=0$ in the original coordinates of (3.6)) then the solution of system (3.10) either satisfies $\rho(t)\equiv 0$ or $\lim_{t\to+\infty}\left(\rho(t)\right)=1$. In this case the solution of (3.6) either tends to an equilibrium point or tends to a periodic solution of period $2\pi/\omega$.

Therefore, we next focus on the case $r(t)>0$ and $\rho(t)>0$. In this case the transformation $r=\exp(x_1)$, $\rho=\exp(x_2)$ with $x_1,x_2\in\mathbb{R}$ allows us to obtain from (3.10) the time-varying system

$$\begin{aligned}&\dot{x}_1(t)=1-\exp\left(2x_1(t)\right)+a\exp\left(x_2(t)\right)\cos\left(\varphi_0+\omega t\right)\\&\dot{x}_2(t)=1-\exp\left(2x_2(t)\right)+b\exp\left(x_1(t)\right)\cos\left(\theta_0+t\right)\\&x(t)=\left(x_1(t),x_2(t)\right)\in\mathbb{R}^2\end{aligned}\tag{3.11}$$

System (3.11) is a system of the form (2.1) with $n=2$ and

$$f(t,x)=\begin{bmatrix}1-\exp\left(2x_1\right)+a\exp\left(x_2\right)\cos\left(\varphi_0+\omega t\right)\\1-\exp\left(2x_2\right)+b\exp\left(x_1\right)\cos\left(\theta_0+t\right)\end{bmatrix}\tag{3.12}$$

When $\omega\in\mathbb{Q}$ the map $f\in C^1\left(\mathbb{R}\times\mathbb{R}^2;\mathbb{R}^2\right)$ defined by (3.12) is periodic with period $T=2m\pi$ for appropriate $m\in\mathbb{N}$. In general, even when $\omega\in\mathbb{R}$ is not rational, the map $f\in C^1\left(\mathbb{R}\times\mathbb{R}^2;\mathbb{R}^2\right)$ defined by (3.12) is almost periodic in $t$ uniformly in $x\in\Lambda$, for every non-empty compact set $\Lambda\subseteq\mathbb{R}^2$.

It is shown in the Appendix that for every $x_0\in\mathbb{R}^2$, $t_0,\theta_0,\varphi_0\in\mathbb{R}$ the solution $x(t)=\left(x_1(t),x_2(t)\right)=\varphi(t,t_0,x_0)$ of (3.11) is defined for all $t\geq t_0$ and satisfies the estimate

$$\limsup_{t\to+\infty}\left(\left|x_i(t)\right|\right)\leq\kappa\left(2\sqrt{\frac{a^2+b^2}{1-|a|-|b|}}\right),\ \text{for } i=1,2,\tag{3.13}$$

where $\kappa(s)=s+\sqrt{2s}$ for $s\geq 0$.

Definition (3.12) gives for all $t\in\mathbb{R}$, $x=\left(x_1,x_2\right)\in\mathbb{R}^2$:

$$\frac{\partial f}{\partial x}(t,x)=\begin{bmatrix}-2\exp\left(2x_1\right) & a\exp\left(x_2\right)\cos\left(\varphi_0+\omega t\right)\\ b\exp\left(x_1\right)\cos\left(\theta_0+t\right) & -2\exp\left(2x_2\right)\end{bmatrix}\tag{3.14}$$

By virtue of (3.7) there exists $B>\kappa\left(2\sqrt{\frac{a^2+b^2}{1-|a|-|b|}}\right)$ so that $\max\left(|a|,|b|\right)\exp\left(B\right)<2$. For the 1-norm of $\mathbb{R}^2$ $|x|=|x|_1=|x_1|+|x_2|$ we get from (3.14) for $x\in C=\left\{\left(z_1,z_2\right)\in\mathbb{R}^2:|z_1|\leq B,|z_2|\leq B\right\}$:

$$
\begin{aligned}
&\mu\left(\frac{\partial f}{\partial x}(t,x)\right) \\
&= \max\left(-2\exp(2x_1)+|b|\exp(x_1)|\cos(\theta_0+t)|, -2\exp(2x_2)+|a|\exp(x_2)|\cos(\varphi_0+\omega t)|\right) \\
&\leq \max\left(-\exp(x_1)\left(2\exp(-B)-|b|\right), -\exp(x_2)\left(2\exp(-B)-|a|\right)\right) \\
&\leq -\exp(-B)\left(2\exp(-B)-\max(|a|,|b|)\right) < 0
\end{aligned} \tag{3.15}
$$

Consequently, since $B > \kappa\left(2\sqrt{\frac{a^2+b^2}{1-|a|-|b|}}\right)$, by virtue of (3.13) and (3.15) all assumptions of Theorem 1 are valid for system (3.11) with arbitrary $\theta_0, \varphi_0 \in \mathbb{R}$, $S(t_0) \equiv \mathbb{R}^2$, $C = \left\{ (z_1, z_2) \in \mathbb{R}^2 : |z_1| \leq B, |z_2| \leq B \right\}$ and $p(t) \equiv \exp(-B)\left(2\exp(-B) - \max(|a|,|b|)\right)$.

Theorem 1 guarantees that for every $\theta_0, \varphi_0 \in \mathbb{R}$ there exists an almost periodic solution $x^* \in C^1(\mathbb{R}; C)$ of (3.11). Moreover, $\lim_{t \to +\infty}\left(\left|\varphi(t,t_0,x_0) - x^*(t)\right|\right) = 0$ for all $t_0 \in \mathbb{R}$, $x_0 \in \mathbb{R}^2$.

Going back to the original coordinates, we have managed to show that for every $z(0) \in \mathbb{R}^4$ with $z_1^2(0) + z_2^2(0) > 0$ and $z_3^2(0) + z_4^2(0) > 0$ there exists an almost periodic solution $z^* \in C^1\left(\mathbb{R}; \mathbb{R}^4\right)$ of (3.6) such that the solution $z(t)$ of (3.6) with initial condition $z(0) \in \mathbb{R}^4$ satisfies $\lim_{t \to +\infty}\left(\left|z(t) - z^*(t)\right|\right) = 0$.

It should be noticed that when $\omega \in \mathbb{Q}$ then Theorem 2 in [15] implies that there exists appropriate $m \in \mathbb{N}$ such that the solution $z^* \in C^1\left(\mathbb{R}; \mathbb{R}^4\right)$ is periodic with period $T = 2m\pi$. $\triangleleft$

## 4. Global Entrainment to Almost Periodic Inputs

In this section we study the global entrainment problem for time-invariant systems of the form

$$
\dot{x} = f(x,u),\ x \in D, u \in \mathbb{R}^m \tag{4.1}
$$

where $D \subseteq \mathbb{R}^n$ is an open set and $f \in C^1\left(D \times \mathbb{R}^m; \mathbb{R}^n\right)$.

We adopt the terminology in [4] and we provide the definition of entrainment to periodic inputs and almost periodic inputs of limited or unlimited magnitude.

**Definition 3:** *Let $D \subseteq \mathbb{R}^n$ be an open set. Consider the system (4.1) with $f \in C^1\left(D \times \mathbb{R}^m; \mathbb{R}^n\right)$. We say that:*

**(a)** *system (4.1) globally entrains to periodic inputs if the following property holds: given an input $u \in C^0\left(\mathbb{R}; \mathbb{R}^m\right)$ which is periodic with period $T > 0$, all solutions of (4.1) converge to a unique periodic solution with period $T > 0$.*

**(b)** *system (4.1) globally entrains to almost periodic inputs if the following property holds: given an input* $u \in C^0(\mathbb{R};\mathbb{R}^m)$ *which is almost periodic, all solutions of (4.1) converge to a unique (almost) periodic solution which has the same basic frequencies with* $u$.

**(c)** *system (4.1) globally entrains to almost periodic inputs of magnitude* $r>0$ *if the following property holds: given an input* $u \in C^0(\mathbb{R};\mathbb{R}^m)$ *with* $\|u\|_\infty \le r$ *which is almost periodic, all solutions of (4.1) converge to a unique almost periodic solution which has the same basic frequencies with* $u$.

**(d)** *system (4.1) globally entrains to periodic inputs of magnitude* $r>0$ *if the following property holds: given an input* $u \in C^0(\mathbb{R};\mathbb{R}^m)$ *with* $\|u\|_\infty \le r$ *which is periodic with period* $T>0$, *all solutions of (4.1) converge to a unique periodic solution with period* $T>0$.

**(e)** *system (4.1) globally entrains to small almost periodic inputs if there exists* $r>0$ *such that system (4.1) globally entrains to almost periodic inputs of magnitude* $r>0$.

**(f)** *system (4.1) globally entrains to small periodic inputs if there exists* $r>0$ *such that system (4.1) globally entrains to periodic inputs of magnitude* $r>0$.

Theorem 1 can be directly applied to the entrainment problem. Applying Theorem 1 we obtain the following corollary.

**Corollary 1:** *Consider system (4.1) with* $f \in C^1(\mathbb{R}^n \times \mathbb{R}^m;\mathbb{R}^n)$ *and* $D=\mathbb{R}^n$. *Let* $|\cdot|$ *denote the standard Euclidean norm of* $\mathbb{R}^n$ *or* $\mathbb{R}^m$. *Suppose that there exist constants* $R,r>0$ *such that the following properties hold:*

**(P1)** *For every* $u \in C^0(\mathbb{R};\mathbb{R}^m)$ *with* $\|u\|_\infty \le r$ *and for every* $x_0 \in \mathbb{R}^n$ *the solution* $x(t)$ *of (4.1) with initial condition* $x(0)=x_0$ *exists for all* $t \ge 0$ *and satisfies* $\limsup_{t\to+\infty}\left(|x(t)|\right) < R$.

**(P2)** *There exists a symmetric, positive definite matrix* $P \in \mathbb{R}^{n\times n}$ *such that the matrix* $Q(x,u) = P\frac{\partial f}{\partial x}(x,u) + \left(\frac{\partial f}{\partial x}(x,u)\right)^T P$ *is negative definite for all* $(x,u) \in \mathbb{R}^n \times \mathbb{R}^m$ *with* $|x| \le R$ *and* $|u| \le r$.

*Then system (4.1) globally entrains to almost periodic inputs of magnitude* $r$. *Moreover, system (4.1) globally entrains to periodic inputs of magnitude* $r$.

The proof of Corollary 1 is provided in Section 6.

**Remarks on Corollary 1: (i)** The importance of Corollary 1 compared with other existence results that are available in the literature (see for example [10] and [5]) is the quantitative and global nature of Corollary 1. Most of the available results are qualitative local results.
**(ii)** Assumption (P1) can be shown by using a radially unbounded function $V \in C^1(\mathbb{R}^n)$ that satisfies the following properties for certain constant $B \in \mathbb{R}$:

$$V(x)\geq B, |u|\leq r \Rightarrow \nabla V(x) f(x,u)<0 \tag{4.2}$$

$$\max\left\{ |x| : V(x)\leq B \right\} < R \tag{4.3}$$

**(iii)** Assumption (P2) involves a localized version of a Demidovich condition which was used by Boris Pavlovich Demidovich for the study of global asymptotic stability of a global solution (see [22]). Assumption (P2) is also met in the study of contraction (see [18]). However, here we require the matrix $Q(x,u)=P\frac{\partial f}{\partial x}(x,u)+\left(\frac{\partial f}{\partial x}(x,u)\right)^T P$ to be negative definite only for an appropriate compact set in $\mathbb{R}^n\times\mathbb{R}^m$ and not globally as in [22] and [18]. Here we avoid the requirement that the Demidovich condition holds globally (which is a very demanding requirement) by employing the ultimate boundedness assumption (P1). In other words, contraction is required only where the trajectories ultimately reside, with ultimate boundedness supplied separately. This feature provides flexibility and can allow the application to cases where a global Demidovich condition cannot be satisfied.

In particular, Corollary 1 implies an interesting qualitative result: see Corollary 2 below. To this end, we use the notion of Input-to-State Stability (ISS) with respect to small inputs, i.e., the property for which $f(0,0)=0$ and there exist $\beta\in KL$, $\gamma\in K$ and a constant $\sigma>0$ such that for every $x_0\in\mathbb{R}^n$ and for every measurable and essentially bounded input $u:\mathbb{R}_+\to\mathbb{R}^m$ with $\|u\|_\infty\leq\sigma$ the solution $x(t)$ of (4.1) with initial condition $x(0)=x_0$ satisfies the following estimate for all $t\geq 0$:

$$|x(t)|\leq \beta\left( |x_0|, t\right)+\gamma\left(\|u\|_\infty\right) \tag{4.4}$$

The Input-to-State Stability (ISS) property with respect to small inputs is studied in the book [21].

**Corollary 2:** *Consider system (4.1) with $f\in C^1\left(\mathbb{R}^n\times\mathbb{R}^m;\mathbb{R}^n\right)$ and $D=\mathbb{R}^n$. Suppose that $f(0,0)=0$ and that system (4.1) satisfies the Input-to-State Stability (ISS) property with respect to small inputs. Moreover, suppose that $0\in\mathbb{R}^n$ is locally exponentially stable (LES) for the input-free system $\dot{x}=f(x,0)$. Then system (4.1) globally entrains to small almost periodic inputs. Moreover, system (4.1) globally entrains to small periodic inputs.*

Corollary 2 is a direct consequence of Corollary 1 and the fact that estimate (4.4) implies the asymptotic estimate $\limsup_{t\to+\infty}\left(|x(t)|\right)\leq\gamma\left(\|u\|_\infty\right)$ for all $u\in C^0\left(\mathbb{R};\mathbb{R}^m\right)$ with $\|u\|_\infty\leq\sigma$. Moreover, since $0\in\mathbb{R}^n$ is locally exponentially stable (LES) for the input-free system $\dot{x}=f(x,0)$, it follows that there exists a symmetric, positive definite matrix $P\in\mathbb{R}^{n\times n}$ such that the matrix $Q(0,0)=P\frac{\partial f}{\partial x}(0,0)+\left(\frac{\partial f}{\partial x}(0,0)\right)^T P$ is negative definite. By continuity of $Q(x,u) = P\frac{\partial f}{\partial x}(x,u) + \left(\frac{\partial f}{\partial x}(x,u)\right)^T P$ it follows that there exists $\delta > 0$ such that the matrix $Q(x,u) = P\frac{\partial f}{\partial x}(x,u) + \left(\frac{\partial f}{\partial x}(x,u)\right)^T P$ is negative definite for all $(x,u) \in \mathbb{R}^n \times \mathbb{R}^m$ with $|x| \le \delta$ and $|u| \le \delta$. It follows that properties (P1) and (P2) of Corollary 1 hold for $R = \delta$ and $r = \min\left(\sigma, \tilde{\delta}, \gamma^{-1}\left(\tilde{\delta}/2\right)\right)$ with $\tilde{\delta} = \min\left(\delta, \sup_{s \ge 0}\left(\gamma(s)\right)\right)$.

## 5. Global Entrainment in Lotka-Volterra Models

Lotka-Volterra systems (see [12]) are dynamical system of the form

$$\dot{y} = diag(y)\left(k + By\right), \tag{5.1}$$

where $B \in \mathbb{R}^{n \times n}$ is a constant real matrix, $k \in \mathbb{R}^n$ is a constant vector and $y \in \operatorname{int}\left(\mathbb{R}^n_+\right)$ is the state. A point $y^* \in \operatorname{int}\left(\mathbb{R}^n_+\right)$ is called an interior (or coexistence) equilibrium point for (5.1) if $k + By^* = 0$. Here we study Lotka-Volterra systems with inputs of the form

$$\dot{y} = diag(y)\left(k + By + Ku\right), \tag{5.2}$$

where $u \in \mathbb{R}^m$ is an input and $K \in \mathbb{R}^{n \times m}$ is a constant matrix. The input $u \in \mathbb{R}^m$ can be used to express the effect of seasonal hunting or food quantity variations on the ecosystem; see for instance [6] and references therein.

When there is an interior equilibrium point $y^*$, the transformation $y_i = diag\left(y_i^*\right)\exp\left(x_i\right)$ for $i = 1,...,n$ brings the Lotka-Volterra system (5.2) to the form

$$\dot{x} = A\left(\exp_n\left(x\right) - 1_n\right) + Ku, \quad x \in \mathbb{R}^n, u \in \mathbb{R}^m, \tag{5.3}$$

where $1_n = \left(1,...,1\right)^T$, $\exp_n(x) = \left(\exp(x_1),...,\exp(x_n)\right)^T$ and $A = Bdiag\left(y^*\right)$. We next focus on system (5.3) that has the equilibrium point $x^* = 0 \in \mathbb{R}^n$ in the input-free case $u \equiv 0$.

The equilibrium point $x^* = 0 \in \mathbb{R}^n$ in the input-free case $u \equiv 0$ is Globally Asymptotically Stable and Locally Exponentially Stable if there exists $p \in \operatorname{int}\left(\mathbb{R}^n_+\right)$ such that the matrix $M = -\frac{1}{2}diag(p)A - \frac{1}{2}A^T diag(p)$ is positive definite (see [12] and [14]). In such a case the matrix $A \in \mathbb{R}^{n \times n}$ is called a Volterra-Lyapunov stable matrix and Global Asymptotic Stability can be proved by exploiting the Volterra-Lyapunov function

$$V(x) = p^T\left(\exp_n(x) - x - 1_n\right) \tag{5.4}$$

Exploiting the function $V(x)$ defined by (5.4) and Corollary 2 we can show the following result.

**Theorem 3:** *Consider the Lotka-Volterra system (5.2) and suppose that there exists an interior equilibrium point* $y^* \in \operatorname{int}\left(\mathbb{R}_+^n\right)$ *with* $k + By^* = 0$. *Moreover, suppose that the matrix* $A = Bdiag\left(y^*\right)$ *is a Volterra-Lyapunov stable matrix, i.e., there exists* $p \in \operatorname{int}\left(\mathbb{R}_+^n\right)$ *such that the matrix* $diag(p)A + A^T diag(p)$ *is negative definite. Then system (5.2) globally entrains to small almost periodic inputs. Moreover, system (5.2) globally entrains to small periodic inputs. Finally, for all almost periodic inputs of sufficiently small magnitude the following equation holds:*

$$\lim_{T \to +\infty} \left( \frac{1}{T} \int_0^T y(t)dt \right) = y^* - B^{-1}K \lim_{T \to +\infty} \left( \frac{1}{T} \int_0^T u(t)dt \right) \tag{5.5}$$

We are not aware of a result similar to that of Theorem 3 in the literature. The biological significance of Theorem 3 is important. It means that ecosystems with a Volterra-Lyapunov stable interaction matrix entrain to almost periodic solutions that are induced by the superposition of seasonal effects and many other periodic effects (e.g., the motion of the moon or the day-night cycle).

The proof of Theorem 3 can be found in Section 6 and provides a methodology of estimating how large almost periodic inputs can be entrained. Moreover, equation (5.5) can be useful for many purposes: equation (5.5) shows how we can affect the mean value of the populations, i.e., $\lim_{T \to +\infty} \left( \frac{1}{T} \int_0^T y(t)dt \right)$, by applying standard techniques (e.g., periodic harvesting).

## 6. Proofs

In this section we provide the proofs of all results that are given above.

**Proof of Theorem 2:** Since $\frac{d}{dt}\left(x^*(t+\tau) - x^*(t)\right) = f\left(t+\tau, x^*(t+\tau)\right) - f\left(t, x^*(t)\right)$ (recall (2.1)), we get for all $t, \tau \in \mathbb{R}$ and $h > 0$:

$$\begin{aligned} &\left| x^*(t+h+\tau) - x^*(t+h) \right| \\ &= \left| x^*(t+\tau) - x^*(t) + h\left( f\left(t+\tau, x^*(t+\tau)\right) - f\left(t, x^*(t+\tau)\right) + f\left(t, x^*(t+\tau)\right) - f\left(t, x^*(t)\right) \right) \right| + o(h) \end{aligned} \tag{6.1}$$

where $\lim_{h \to 0^+} \left( \frac{o(h)}{h} \right) = 0$. Since $x^*(t) \in C$ for all $t \in \mathbb{R}$, we obtain from (6.1):

$$\begin{aligned} &\left| x^*(t+h+\tau) - x^*(t+h) \right| \\ &\leq \left| x^*(t+\tau) - x^*(t) + h\left( f\left(t, x^*(t+\tau)\right) - f\left(t, x^*(t)\right) \right) \right| + h\left| f\left(t+\tau, x^*(t+\tau)\right) - f\left(t, x^*(t+\tau)\right) \right| + o(h) \\ &\leq \left| x^*(t+\tau) - x^*(t) + h\left( \int_0^1 \frac{\partial f}{\partial x}\left(t, x^*(t) + \lambda\left(x^*(t+\tau) - x^*(t)\right)\right) d\lambda \right)\left(x^*(t+\tau) - x^*(t)\right) \right| + hF(\tau) + o(h) \end{aligned} \tag{6.2}$$

where

$$F(\tau) := \sup\left\{ \left| f(l+\tau,z) - f(l,z) \right| : l \in \mathbb{R}, z \in C \right\} \tag{6.3}$$

Thus, we get from (6.2) for all $t,\tau \in \mathbb{R}$ and $h > 0$:

$$\begin{aligned} &\left| x^*(t+h+\tau) - x^*(t+h) \right| \\ &\leq \left\| I_n + h\left( \int_0^1 \frac{\partial f}{\partial x}\left(t, x^*(t) + \lambda\left(x^*(t+\tau) - x^*(t)\right)\right) d\lambda \right) \right\| \left| x^*(t+\tau) - x^*(t) \right| + hF(\tau) + o(h) \end{aligned} \tag{6.4}$$

where $\|\cdot\|$ is the induced matrix norm. Using (6.4) and definition (1.1) we get for all $t,\tau \in \mathbb{R}$:

$$\begin{aligned} &\limsup_{h\to 0^+} \left( \frac{\left| x^*(t+h+\tau) - x^*(t+h) \right| - \left| x^*(t+\tau) - x^*(t) \right|}{h} \right) \\ &\leq \mu\left( \int_0^1 \frac{\partial f}{\partial x}\left(t, x^*(t) + \lambda\left(x^*(t+\tau) - x^*(t)\right)\right) d\lambda \right) \left| x^*(t+\tau) - x^*(t) \right| + F(\tau) \end{aligned} \tag{6.5}$$

Since $x^*(t) \in C$ for all $t \in \mathbb{R}$, due to convexity of the set $C \subseteq \mathbb{R}^n$ it holds that $\left( x^*(t) + \lambda\left(x^*(t+\tau) - x^*(t)\right) \right) \in C$ for all $\lambda \in [0,1]$. Employing Lemma 2 in [15] and (2.2) we get from (6.5) for all $t,\tau \in \mathbb{R}$:

$$\limsup_{h\to 0^+} \left( \frac{\left| x^*(t+h+\tau) - x^*(t+h) \right| - \left| x^*(t+\tau) - x^*(t) \right|}{h} \right) \leq -p(t)\left| x^*(t+\tau) - x^*(t) \right| + F(\tau) \tag{6.6}$$

Exploiting Lemma 2.5 (Comparison Lemma) in [16], we get from (6.6) for all $t,\tau,t_0 \in \mathbb{R}$ with $t \geq t_0$:

$$\left| x^*(t+\tau) - x^*(t) \right| \leq \left| x^*(t_0+\tau) - x^*(t_0) \right| \exp\left( -\int_{t_0}^{t} p(s)ds \right) + F(\tau) \int_{t_0}^{t} \exp\left( -\int_{l}^{t} p(s)ds \right) dl \tag{6.7}$$

Define for all $t \in \mathbb{R}$:

$$y(t) = x^*(-t) \tag{6.8}$$

Since $\frac{d}{ds}\left( y(s) - y(s-\tau) \right) = -f\left(-s, y(s)\right) + f\left(-s+\tau, y(s-\tau)\right)$ (recall (2.1) and (6.8)), we get for all $s,\tau \in \mathbb{R}$ and $h > 0$:

$$\begin{aligned} &\left| y(s+h) - y(s-\tau+h) \right| \\ &= \left| y(s) - y(s-\tau) + h\left( -f\left(-s, y(s)\right) + f\left(-s, y(s-\tau)\right) - f\left(-s, y(s-\tau)\right) + f\left(-s+\tau, y(s-\tau)\right) \right) \right| + o(h) \end{aligned} \tag{6.9}$$

where $\lim_{h\to 0^+}\left(\frac{o(h)}{h}\right)=0$. Since $y(t)\in C$ for all $t\in\mathbb{R}$ (recall (6.8) and the fact that $x^*(t)\in C$) for all $t\in\mathbb{R}$, we obtain from (6.9) and definition (6.3):

$$\begin{aligned}&|y(s+h)-y(s-\tau+h)|\\&\le\left|y(s)-y(s-\tau)+h\left(-f\left(-s,y(s)\right)+f\left(-s,y(s-\tau)\right)\right)\right|\\&+h\left|f\left(-s+\tau,y(s-\tau)\right)-f\left(-s,y(s-\tau)\right)\right|+o(h)\\&\le\left|y(s)-y(s-\tau)+h\left(-\int_0^1\frac{\partial f}{\partial x}\left(-s,y(s-\tau)+\lambda\left(y(s)-y(s-\tau)\right)\right)d\lambda\right)\left(y(s)-y(s-\tau)\right)\right|\\&+hF(\tau)+o(h)\end{aligned}\tag{6.10}$$

Thus, we get from (6.10) for all $s,\tau\in\mathbb{R}$ and $h>0$:

$$\begin{aligned}&|y(s+h)-y(s-\tau+h)|\\&\le\left\|I_n+h\left(-\int_0^1\frac{\partial f}{\partial x}\left(-s,y(s-\tau)+\lambda\left(y(s)-y(s-\tau)\right)\right)d\lambda\right)\right\|\left|y(s)-y(s-\tau)\right|+hF(\tau)+o(h)\end{aligned}\tag{6.11}$$

where $\|\cdot\|$ is the induced matrix norm. Using (6.11) and definition (1.1) we get for all $s,\tau\in\mathbb{R}$:

$$\begin{aligned}&\limsup_{h\to 0^+}\left(\frac{|y(s+h)-y(s-\tau+h)|-|y(s)-y(s-\tau)|}{h}\right)\\&\le\mu\left(-\int_0^1\frac{\partial f}{\partial x}\left(-s,y(s-\tau)+\lambda\left(y(s)-y(s-\tau)\right)\right)d\lambda\right)\left|y(s)-y(s-\tau)\right|+F(\tau)\end{aligned}\tag{6.12}$$

Since $y(t)\in C$ for all $t\in\mathbb{R}$ (recall (6.8) and the fact that $x^*(t)\in C$ for all $t\in\mathbb{R}_-$, due to convexity of the set $C\subseteq\mathbb{R}^n$ it holds that $\left(y(s-\tau)+\lambda\left(y(s)-y(s-\tau)\right)\right)\in C$ for all $\lambda\in[0,1]$. Employing Lemma 2 in [15] and (2.2) we get from (6.12) for all $s,\tau\in\mathbb{R}$:

$$\limsup_{h\to 0^+}\left(\frac{|y(s+h)-y(s-\tau+h)|-|y(s)-y(s-\tau)|}{h}\right)\le-\beta(-s)\left|y(s)-y(s-\tau)\right|+F(\tau)\tag{6.13}$$

Exploiting once again Lemma 2.5 (Comparison Lemma) in [16] we get from (6.13) for all $s,\tau,s_0\in\mathbb{R}$ with $s\ge s_0$:

$$|y(s)-y(s-\tau)|\le|y(s_0)-y(s_0-\tau)|\exp\left(-\int_{s_0}^{s}\beta(-r)dr\right)+F(\tau)\int_{s_0}^{s}\exp\left(-\int_{l}^{s}\beta(-r)dr\right)dl\tag{6.14}$$

Consequently, by setting $s=-t$, $s_0=-t_0$, we get from (6.14) and definition (6.8) for all $t,\tau,t_0\in\mathbb{R}$ with $t\le t_0$:

$$\left|x^*(t+\tau)-x^*(t)\right| \le \left|x^*(t_0+\tau)-x^*(t_0)\right| \exp\left(-\int_t^{t_0}\beta(r)dr\right)+F(\tau)\int_t^{t_0}\exp\left(-\int_t^{l}\beta(r)dr\right)dl \tag{6.15}$$

Let arbitrary $t\in\mathbb{R}$ be given. If property (I) of Theorem 1 holds then estimate (6.7) implies that there exists $t_0\in(-\infty,t)$ with

$$\begin{aligned}&\left|x^*(t+\tau)-x^*(t)\right| \le \left|x^*(t_0+\tau)-x^*(t_0)\right|\exp(-c)+KF(\tau)\\ &\le \left\|\delta_\tau x^*-x^*\right\|_\infty \exp(-c)+KF(\tau)\end{aligned} \tag{6.16}$$

If property (II) holds then estimate (6.15) implies that there exists $t_0\in(t,+\infty)$ for which (6.16) holds. Consequently, we obtain for all $t,\tau\in\mathbb{R}$:

$$\left|x^*(t+\tau)-x^*(t)\right| \le \left\|\delta_\tau x^*-x^*\right\|_\infty \exp(-c)+KF(\tau) \tag{6.17}$$

Estimate (6.17) implies the following estimate for all $\tau\in\mathbb{R}$:

$$\left\|\delta_\tau x^*-x^*\right\|_\infty \le \frac{K}{1-\exp(-c)}F(\tau) \tag{6.18}$$

The fact that $f:\mathbb{R}\times\mathbb{R}^n\to\mathbb{R}^n$ is almost periodic in $t$ uniformly in $x\in C$ in conjunction with Fact I imply that $f$ has a set $\Omega\subset\mathbb{R}$ of basic frequencies, i.e., $\Omega\subset\mathbb{R}$ is a countable set which is linearly independent over the rationals and for every sequence $\{\tau_k:k=1,2,...\}\subseteq\mathbb{R}$ with $\lim_{k\to+\infty}\left(\exp\left(i\omega_p\tau_k\right)\right)=1$ for each $\omega_p\in\Omega$ it holds that $\lim_{k\to+\infty}\left(\sup_{t\in\mathbb{R},x\in C}\left(\left|f(t+\tau_k,x)-f(t,x)\right|\right)\right)=0$. Using definition (6.3), we conclude that for every sequence $\{\tau_k:k=1,2,...\}\subseteq\mathbb{R}$ with $\lim_{k\to+\infty}\left(\exp\left(i\omega_p\tau_k\right)\right)=1$ for each $\omega_p\in\Omega$ it holds that $\lim_{k\to+\infty}\left(F(\tau_k)\right)=0$.

Therefore, estimate (6.18) implies that for every sequence $\{\tau_k:k=1,2,...\}\subseteq\mathbb{R}$ with $\lim_{k\to+\infty}\left(\exp\left(i\omega_p\tau_k\right)\right)=1$ for each $\omega_p\in\Omega$ it holds that $\lim_{k\to+\infty}\left(\left\|\delta_{\tau_k}x^*-x^*\right\|_\infty\right)=\lim_{k\to+\infty}\left(\sup_{t\in\mathbb{R}}\left(\left|x^*(t+\tau_k)-x^*(t)\right|\right)\right)=0$. Fact II implies that $x^*$ is almost periodic and $\Omega\subset\mathbb{R}$ is the set of basic frequencies of $x^*$. Consequently, $x^*$ has the same basic frequencies with $f$. The proof is complete. $\triangleleft$

**Proof of Lemma 1:** By virtue of the fact that $f\in C^1\left(\mathbb{R}\times\mathbb{R}^n;\mathbb{R}^n\right)$ is almost periodic in $t$ uniformly in $x\in C$, there exists a sequence $\{\tau_m\ge 0:m=1,...\}$ with $\lim_{m\to+\infty}(\tau_m)=+\infty$ and $\lim_{m\to+\infty}\left(\sup\left\{\left|f(t+\tau_m,x)-f(t,x)\right|:t\in\mathbb{R},x\in C\right\}\right)=0$.

Let arbitrary $T>0$ be given. Since $\lim_{m\to+\infty}(\tau_m)=+\infty$ there exists $\bar{m}\in\mathbb{N}$ such that $\tau_m\ge T+\max(t_0,\tau)$ for all $m\ge\bar{m}$. Consequently, $\varphi(t+\tau_m,t_0,x_0)\in C$ for all $t\ge -T$ and $m\ge\bar{m}$.

Define for all $m\ge\bar{m}$ the sequence of functions $y_m:[t_0-T,t_0+T]\to C$ by means of the formula

$$y_m(t) = \varphi\left(t + \tau_m, t_0, x_0\right) \tag{6.19}$$

Since $\dot{y}_m(t) = f\left(t + \tau_m, y_m(t)\right)$ for all $m \geq \bar{m}$, $t \in \left[t_0 - T, t_0 + T\right]$ and since $\sup\left\{ \left| f(t,x)\right| : t \in \mathbb{R}, x \in C\right\} < +\infty$ (recall assumption (A)) we obtain the existence of a constant $M > 0$ such that $\sup_{t\in[t_0-T,t_0+T]}\left(\left|\dot{y}_m(t)\right|\right) \leq M$ for all $m \geq \bar{m}$.

Therefore, the sequence of continuous functions $\left\{y_m \in C^0\left(\left[t_0 - T, t_0 + T\right]; C\right) : m \geq \bar{m}\right\}$ is uniformly bounded (since $C$ is compact) and equicontinuous (since $\sup_{t\in[t_0-T,t_0+T]}\left(\left|\dot{y}_m(t)\right|\right) \leq M$ ). By virtue of the Arzela-Ascoli theorem there exists $x_T \in C^0\left(\left[t_0 - T, t_0 + T\right]; \mathbb{R}^n\right)$ and a subsequence $\left\{m_k \geq \bar{m} : k = 1,...\right\}$ with $\sup_{t\in[t_0-T,t_0+T]}\left(\left|y_{m_k}(t) - x_T(t)\right|\right) \to 0$. Since $C \subseteq \mathbb{R}^n$ is closed it follows that $x_T \in C^0\left(\left[t_0 - T, t_0 + T\right]; C\right)$.

Since $\dot{y}_m(t) = f\left(t + \tau_m, y_m(t)\right)$ for all $m \geq \bar{m}$, $t \in \left[t_0 - T, t_0 + T\right]$, we also have for all $m \geq \bar{m}$ and $t \in \left[t_0 - T, t_0 + T\right]$:

$$\begin{aligned}
&\left| x_T(t) - x_T(t_0) - \int_{t_0}^{t} f\left(s, x_T(s)\right) ds \right| \\
&= \left| x_T(t) - y_m(t) - x_T(t_0) + y_m(t_0) - \int_{t_0}^{t}\left( f\left(s, x_T(s)\right) - f\left(s + \tau_m, y_m(s)\right)\right) ds \right| \\
&\leq \left| x_T(t) - y_m(t)\right| + \left| x_T(t_0) - y_m(t_0)\right| + \int_{\min(t,t_0)}^{\max(t,t_0)} \left| f\left(s, x_T(s)\right) - f\left(s + \tau_m, y_m(s)\right)\right| ds
\end{aligned}$$

Assumption (A) guarantees the existence of a constant $L > 0$ that satisfies $\left| f(t,x) - f(t,y)\right| \leq L\left|x - y\right|$ for all $t \in \mathbb{R}$, $x, y \in C$. Thus, we get for all $m \geq \bar{m}$, $t \in \left[t_0 - T, t_0 + T\right]$:

$$\begin{aligned}
&\left| x_T(t) - x_T(t_0) - \int_{t_0}^{t} f\left(s, x_T(s)\right) ds \right| \leq 2 \sup_{s\in[t_0-T,t_0+T]}\left(\left| y_m(s) - x_T(s)\right|\right) \\
&+ \int_{\min(t,t_0)}^{\max(t,t_0)} \left| f\left(s, x_T(s)\right) - f\left(s + \tau_m, x_T(s)\right) + f\left(s + \tau_m, x_T(s)\right) - f\left(s + \tau_m, y_m(s)\right)\right| ds \\
&\leq 2 \sup_{s\in[t_0-T,t_0+T]}\left(\left| y_m(s) - x_T(s)\right|\right) + \int_{\min(t,t_0)}^{\max(t,t_0)} \left| f\left(s, x_T(s)\right) - f\left(s + \tau_m, x_T(s)\right)\right| ds \\
&+ \int_{\min(t,t_0)}^{\max(t,t_0)} \left| f\left(s + \tau_m, x_T(s)\right) - f\left(s + \tau_m, y_m(s)\right)\right| ds \\
&\leq \left(2 + L\left|t - t_0\right|\right) \sup_{s\in[t_0-T,t_0+T]}\left(\left| y_m(s) - x_T(s)\right|\right) + \int_{\min(t,t_0)}^{\max(t,t_0)} \left| f\left(s, x_T(s)\right) - f\left(s + \tau_m, x_T(s)\right)\right| ds \\
&\leq \left(2 + L\left|t - t_0\right|\right) \sup_{s\in[t_0-T,t_0+T]}\left(\left| y_m(s) - x_T(s)\right|\right) + \left|t - t_0\right| \sup\left\{ \left| f(l + \tau_m, x) - f(l, x)\right| : l \in \mathbb{R}, x \in C\right\}
\end{aligned}$$

Therefore, we get for all $m \geq \bar{m}$, $t \in [t_0 - T, t_0 + T]$:

$$\sup_{t\in[t_0-T,t_0+T]}\left(\left|x_T(t)-x_T(t_0)-\int_{t_0}^{t} f(s,x_T(s))ds\right|\right)$$
$$\leq (2+LT)\sup_{s\in[t_0-T,t_0+T]}\left(\left|y_m(s)-x_T(s)\right|\right)+T\sup\left\{\,\left|f(l+\tau_m,x)-f(l,x)\right|:l\in\mathbb{R},x\in C\right\}$$

Exploiting the subsequence $\{m_k \geq \bar{m} : k = 1,...\}$ with $\sup_{t\in[t_0-T,t_0+T]}\left(\left|y_{m_k}(t)-x_T(t)\right|\right)\to 0$ and the fact that $\lim_{m\to+\infty}\left(\sup\left\{\,\left|f(t+\tau_m,x)-f(t,x)\right|:t\in\mathbb{R},x\in C\right\}\right)=0$, we obtain that $x_T(t)=x_T(t_0)+\int_{t_0}^{t} f(s,x_T(s))ds$ for all $t\in[t_0-T,t_0+T]$. Consequently, $x_T\in C^1([t_0-T,t_0+T];C)$ is a solution of (2.1) and $x_T(t)=\varphi(t,t_0,x_T(t_0))$ for $t\in[t_0-T,t_0+T]$.

Consider the sequence

$$\{x_T(t_0)\in C : T=1,2,...\}\subseteq C \tag{6.20}$$

The sequence $\{x_T(t_0)\in C : T=1,2,...\}$ is bounded. By virtue of the Bolzano–Weierstrass theorem and compactness of $C$, there exists $\bar{x}\in C$ and a subsequence $\{T_k\in\mathbb{N} : k=1,...\}$ with $T_k\to+\infty$ and $x_{T_k}(t_0)\to\bar{x}$.

We show next that $\varphi(t,t_0,\bar{x})$ is defined for all $t\in\mathbb{R}$ and satisfies $\varphi(t,t_0,\bar{x})\in C$ for all $t\in\mathbb{R}$.

Clearly, there exists $t_{\max}\in(t_0,+\infty]$ and $t_{\min}\in[-\infty,t_0)$ for which $\varphi(t,t_0,\bar{x})$ is defined for all $t\in(t_{\min},t_{\max})$. Moreover, if $t_{\max}<+\infty$ then $\limsup_{t\to t_{\max}^-}\left(\left|\varphi(t,t_0,\bar{x})\right|\right)=+\infty$. Furthermore, if $t_{\min}>-\infty$ then $\limsup_{t\to t_{\min}^+}\left(\left|\varphi(t,t_0,\bar{x})\right|\right)=+\infty$.

Let arbitrary $\tau\in[t_0,t_{\max})$ be given. Let $\bar{k}\in\mathbb{N}$ be such that $t_0+T_k>\tau$ for $k\geq\bar{k}$. It follows that $x_{T_k}(t)=\varphi(t,t_0,x_{T_k}(t_0))$ for all $t\in[t_0,\tau]$ and $k\geq\bar{k}$. Since $x_{T_k}(t_0)\to\bar{x}$, by continuity of the solution with respect to the initial conditions we have $x_{T_k}(\tau)=\varphi(\tau,t_0,x_{T_k}(t_0))\to\varphi(\tau,t_0,\bar{x})$. Since $x_{T_k}(\tau)\in C$ and since $C\subseteq\mathbb{R}^n$ is closed it follows that $\varphi(\tau,t_0,\bar{x})\in C$. Since $\tau\in[t_0,t_{\max})$ is arbitrary we get $\varphi(t,t_0,\bar{x})\in C$ for all $t\in[t_0,t_{\max})$.

Therefore, a contradiction argument and the fact that $C$ is bounded guarantees that $t_{\max}=+\infty$.

A similar argument shows that $t_{\min}=-\infty$ and $\varphi(t,t_0,\bar{x})\in C$ for all $t\leq t_0$. The conclusion of the lemma holds with $x^*(t)=\varphi(t,t_0,\bar{x})$ for $t\in\mathbb{R}$. The proof is complete. $\triangleleft$

**Proof of Lemma 2:** Let arbitrary $x_0 \in S(t_0)$ be given. Using (2.1) we obtain for all $t \geq t_0$:

$$\begin{aligned}&\frac{\partial}{\partial t}\left(\varphi(t,t_0,x_0)-x^*(t)\right)=f\left(t,\varphi(t,t_0,x_0)\right)-f\left(t,x^*(t)\right)\\&=\left(\int_0^1\frac{\partial f}{\partial x}\left(t,(1-\lambda)x^*(t)+\lambda\varphi(t,t_0,x_0)\right)d\lambda\right)\left(\varphi(t,t_0,x_0)-x^*(t)\right)\end{aligned} \tag{6.21}$$

Exploiting the Coppel inequality (see [24]) and (6.21) we obtain for all $t \geq T = \max(\tau, t_0)$:

$$\begin{aligned}&\left|\varphi(t,t_0,x_0)-x^*(t)\right|\\&\leq\left|\varphi(T,t_0,x_0)-x^*(T)\right|\exp\left(\int_T^t\mu\left(\int_0^1\frac{\partial f}{\partial x}\left(s,(1-\lambda)x^*(s)+\lambda\varphi(s,t_0,x_0)\right)d\lambda\right)ds\right)\end{aligned} \tag{6.22}$$

Since $\varphi(t,t_0,x_0)\in C$ and $x^*(t)\in C$ for all $t\geq\tau$, due to convexity of the set $C\subseteq\mathbb{R}^n$ it holds that $\left((1-\lambda)x^*(s)+\lambda\varphi(s,t_0,x_0)\right)\in C$ for all $\lambda\in[0,1]$ and $s\geq T$. Since $\mu\left(\frac{\partial f}{\partial z}(t,z)\right)\leq -p(t)$ for all $t\geq 0$, $z\in C$, it follows that $\mu\left(\frac{\partial f}{\partial x}\left(s,(1-\lambda)x^*(s)+\lambda\varphi(s,t_0,x_0)\right)\right)\leq -p(s)$ for all $\lambda\in[0,1]$ and $s\geq T$. Exploiting Lemma 2 in [15] and (6.22) we get for all $t\geq T$:

$$\left|\varphi(t,t_0,x_0)-x^*(t)\right|\leq\left|\varphi(T,t_0,x_0)-x^*(T)\right|\exp\left(-\int_T^t p(s)\,ds\right) \tag{6.23}$$

Estimate (6.23) and the fact $\int_0^{+\infty}p(s)\,ds=+\infty$ imply that $\lim_{t\to+\infty}\left(\left|\varphi(t,t_0,x_0)-x^*(t)\right|\right)=0$. The proof is complete. ◁

**Proof of Proposition 1:** The fact that $\int_0^{+\infty}p(s)\,ds=+\infty$ is a consequence of the fact that $\int_0^{+\infty}p(s)\,ds=\sum_{k=0}^{+\infty}\int_{kT}^{(k+1)T}p(s)\,ds=\sum_{k=0}^{+\infty}\int_0^T p(kT+s)\,ds$. Due to the fact that $p\in C^0(\mathbb{R})$ is a $T-$periodic function we get that $p(kT+s)=p(s)$ for all $k\in\mathbb{Z}$ and $s\in\mathbb{R}$. Consequently, $\int_0^{+\infty}p(s)\,ds=\sum_{k=0}^{+\infty}\int_0^T p(s)\,ds=\int_0^T p(s)\,ds\sum_{k=0}^{+\infty}1=+\infty$.

Define for $t\geq t_0$:

$$g(t,t_0)=\int_{t_0}^t p(s)ds \text{ and } r:=\int_0^T p(s)ds>0 \tag{6.24}$$

For $m=\left\lfloor \frac{t_0}{T} \right\rfloor$ and $k=\left\lfloor \frac{t}{T} \right\rfloor \geq m$ we get from definition (6.24):

$$\begin{aligned} g(t,t_0) &= \int_{mT}^{kT} p(s)ds - \int_{mT}^{t_0} p(s)ds + \int_{kT}^{t} p(s)ds \\ &= \int_0^{(k-m)T} p(s+mT)ds - \int_0^{t_0-mT} p(s+mT)ds + \int_0^{t-kT} p(s+kT)ds \end{aligned} \tag{6.25}$$

Since $p \in C^0(\mathbb{R})$ is a $T-$periodic function we get from (6.25):

$$g(t,t_0) = \int_0^{(k-m)T} p(s)ds - \int_0^{t_0-mT} p(s)ds + \int_0^{t-kT} p(s)ds \tag{6.26}$$

We assume next that $k \geq m+1$. In this case we get from (6.26):

$$\begin{aligned} g(t,t_0) &= \sum_{l=0}^{k-m-1} \int_{lT}^{(l+1)T} p(s)ds - \int_0^{t_0-mT} p(s)ds + \int_0^{t-kT} p(s)ds \\ &= \sum_{l=0}^{k-m-1} \int_0^{T} p(s+lT)ds - \int_0^{t_0-mT} p(s)ds + \int_0^{t-kT} p(s)ds \end{aligned} \tag{6.27}$$

Since $p \in C^0(\mathbb{R})$ is a $T-$periodic function we get from (6.27) and (6.24):

$$\begin{aligned} g(t,t_0) &= \sum_{l=0}^{k-m-1} \int_0^{T} p(s)ds - \int_0^{t_0-mT} p(s)ds + \int_0^{t-kT} p(s)ds \\ &= (k-m)r - \int_0^{t_0-mT} p(s)ds + \int_0^{t-kT} p(s)ds \end{aligned} \tag{6.28}$$

Equation (6.25) implies that (6.28) is valid even if $k=m$. Define:

$$A = \max_{\tau\in[0,T]} \left( \left| \int_0^{\tau} p(s)ds \right| \right) \tag{6.29}$$

It follows from (6.28), (6.29) and the fact that $m=\left\lfloor \frac{t_0}{T} \right\rfloor$ and $k=\left\lfloor \frac{t}{T} \right\rfloor$ (which implies that $\frac{t-t_0}{T} - 1 \leq k-m \leq \frac{t-t_0}{T}+1$) that the following estimates hold:

$$\left( \frac{t-t_0}{T} - 1 \right) r - 2A \leq g(t,t_0) \leq \left( \frac{t-t_0}{T} + 1 \right) r + 2A \tag{6.30}$$

Using (6.24), (6.29) and (6.30) we get:

$$\int_{t_0}^{t} \exp\left(-\int_{l}^{t} p(s)ds\right)dl = \int_{t_0}^{t} \exp(-g(t,l))\,dl \le \exp(2A+r)\int_{t_0}^{t} \exp\left(-\left(\frac{t-l}{T}\right)r\right)dl$$
$$= \exp(2A+r)\frac{1-\exp\left(-\frac{t-t_0}{T}r\right)}{r}T \le \frac{T}{r}\exp(2A+r) \le \frac{T\exp(3A)}{r} \tag{6.31}$$

Consequently, the inequality $\int_{t_0}^{t} \exp\left(-\int_{l}^{t} p(s)ds\right)dl \le K$ holds for all $t \ge t_0$ with $K = \frac{T\exp(3A)}{r}$.

Using (6.30) we conclude that the inequality $\int_{t_0}^{t} p(s)ds \ge c$ holds with $c=1$ for $t_0 \le t - T - \frac{(1+2A)T}{r}$. The proof is complete. ⊲

**Proof of Corollary 1:** We apply Theorem 1 with $C = \left\{ x \in \mathbb{R}^n : |x| \le R \right\}$, which is a convex, compact set.

Let $u \in C^0\left(\mathbb{R};\mathbb{R}^m\right)$ with $\|u\|_\infty \le r$ be a given (arbitrary) almost periodic input. Then $\tilde{f} \in C^1\left(\mathbb{R}\times\mathbb{R}^n;\mathbb{R}^n\right)$ defined by $\tilde{f}(t,x) := f(x,u(t))$ satisfies assumption (A) and is almost periodic in $t$ uniformly in $x \in C$. Moreover, $\tilde{f}$ has the same basic frequencies with $u$.

Let $\varphi(t,t_0,x_0)$ denote the unique solution of $\dot{x}(t) = \tilde{f}(t,x(t)) = f(x(t),u(t))$ with initial condition $x(t_0) = x_0$. Let $y(t)$ be the solution of (4.1) with initial condition $y(0) = x_0$ corresponding to input $\delta_{t_0}u$. It follows that $y(t-t_0) = \varphi(t,t_0,x_0)$ for all $t \ge t_0$. This identification of solutions is being used in what follows.

Let $S(t_0) \subseteq \mathbb{R}^n$ be the set of all $x_0 \in \mathbb{R}^n$ for which there exists $\tau \in \mathbb{R}$ with $\varphi(t,t_0,x_0) \in C$ for all $t \ge \tau$. Property (P1) guarantees that $S(t_0) = \mathbb{R}^n$ for all $t_0 \in \mathbb{R}$.

Let $L \in \mathbb{R}^{n\times n}$ invertible matrix such that $P = L^T L$. For the weighted 2-norm $|x| = |Lx|_2$, where $|\cdot|_2$ denotes the standard Euclidean norm in $\mathbb{R}^n$, we get $\max\left\{ y^T\left(PG + G^T P\right)y : y^T P y = 1\right\} = 2\mu(G)$ for every matrix $G \in \mathbb{R}^{n\times n}$. Since $Q(x,u) = P\frac{\partial f}{\partial x}(x,u) + \left(\frac{\partial f}{\partial x}(x,u)\right)^T P$, we get for all $t \in \mathbb{R}$ and $z \in C$:

$$2\mu\left(\frac{\partial \tilde{f}}{\partial z}(t,z)\right) = 2\mu\left(\frac{\partial f}{\partial z}(z,u(t))\right)$$
$$= \max\left\{ y^T Q(z,u(t))\,y : y^T P y = 1\right\} \le \max\left\{ y^T Q(x,v)\,y : y^T P y = 1, |x| \le R, |v| \le r\right\} = -b \tag{6.32}$$

Assumption (P2) guarantees that $b>0$. It follows from inequality (6.32) that all assumptions of Theorem 3 hold with $p(t)\equiv b/2>0$. More specifically, property (I) of Theorem 1 holds for every $t\in\mathbb{R}$ with $t_0=t-1$, $c=b$ and $K=2/b$.

Consequently, there exists an almost periodic solution $x^*\in C^1(\mathbb{R};C)$ of $\dot{x}(t)=\tilde{f}\left(t,x(t)\right)=f\left(x(t),u(t)\right)$ which has the same basic frequencies with $\tilde{f}$, i.e, the same basic frequencies with $u$. Moreover, $\lim_{t\to+\infty}\left(\left|\varphi\left(t,t_0,x_0\right)-x^*(t)\right|\right)=0$ for all $t_0\in\mathbb{R}$, $x_0\in\mathbb{R}^n$.

The proof for the periodic case is exactly the same but we use Theorem 2 in [15] instead of Theorem 3. The proof is complete. $\triangleleft$

**Proof of Theorem 3:** Due to Corollary 2, it suffices to show that system (5.3) satisfies the ISS property with respect to small inputs.

The derivative of the function $V(x)$ defined by (5.4) along the trajectories of system (5.3) satisfies for all $(x,u)\in\mathbb{R}^n\times\mathbb{R}^m$

$$\begin{aligned}\dot{V}(x,u)&=-\left(\exp_n(x)-1_n\right)^T M\left(\exp_n(x)-1_n\right)+\left(\exp_n(x)-1_n\right)^T diag\left(p\right)Ku\\ &\le -\frac{1}{2}\left(\exp_n(x)-1_n\right)^T M\left(\exp_n(x)-1_n\right)+\frac{1}{2}u^T K^T diag\left(p\right)M^{-1}diag\left(p\right)Ku\\ &\le -c\left|\exp_n(x)-1_n\right|^2+G\left|u\right|^2\end{aligned} \tag{6.33}$$

where $|\cdot|$ denotes the standard Euclidean norm of $\mathbb{R}^n$ or $\mathbb{R}^m$, $c>0$ is a constant for which the inequality $\xi^T M\xi\ge 2c|\xi|^2$ holds for all $\xi\in\mathbb{R}^n$ and $G:=\frac{1}{2}\left\|K^T diag\left(p\right)M^{-1}diag\left(p\right)K\right\|$.

It is clear that the function $V(x)$ defined by (5.4) satisfies the following estimate for all $x\in\mathbb{R}^n$:

$$V(x)\le p^T 1_n \max_{i=1,\dots,n}\left(g(x_i)\right) \tag{6.34}$$

where the function $g:\mathbb{R}\to\mathbb{R}_+$ is the positive definite function defined by the equation

$$g(y):=\exp(y)-y-1 \tag{6.35}$$

The following inequality holds for all $y\in\mathbb{R}$:

$$\left|\exp(y)-1\right|\ge\min\left(g(y),1/2\right) \tag{6.36}$$

The function $\kappa:\mathbb{R}_+\to[0,1)$ is an increasing, continuous function that satisfies $\kappa(0)=0$, $\lim_{s\to+\infty}\left(\kappa\left(s\right)\right)=1$ and $\left|\exp(y)-1\right|\ge\kappa\left(g(y)\right)$ for all $y\in\mathbb{R}$. Combining (6.33), (6.34) and (6.36) we get:

$$\dot{V}(x,u) \le -c\sum_{i=1}^{n}\left|\exp(x_i)-1\right|^2 + G|u|^2 \le -c\sum_{i=1}^{n}\left(\min\left(g(x_i),1/2\right)\right)^2 + G|u|^2$$
$$\le -c\left(\min\left(\max_{i=1,\dots,n}\left(g(x_i)\right),1/2\right)\right)^2 + G|u|^2 \le -c\left(\min\left(\frac{V(x)}{p^T 1_n},1/2\right)\right)^2 + G|u|^2 \tag{6.37}$$

It becomes clear from (6.37) that the inequality $|u| \le \sqrt{\frac{c}{2G}}\min\left(\frac{V(x)}{p^T 1_n},1/2\right)$ implies the inequality $\dot{V}(x,u) \le -\frac{c}{2}\left(\min\left(\frac{V(x)}{p^T 1_n},1/2\right)\right)^2$. Exploiting this fact and Lemma 2.14 on page 82 in [13], we conclude that for every $\sigma < \frac{1}{2}\sqrt{\frac{c}{2G}}$ there exist $\beta \in KL$, $\gamma \in K$ such that for every $x_0 \in \mathbb{R}^n$ and for every measurable and essentially bounded input $u:\mathbb{R}_+ \to \mathbb{R}^m$ with $\|u\|_\infty \le \sigma$ the solution $x(t)$ of (5.3) with initial condition $x(0)=x_0$ satisfies estimate (4.4) for all $t \ge 0$.

Every solution $y(t)$ of (5.2) with $y(0) \in \operatorname{int}\left(\mathbb{R}^n_+\right)$ corresponds through the transformation $y_i = diag\left(y_i^*\right)\exp\left(x_i\right)$ for $i=1,\dots,n$ to a solution of (5.3). Since for every measurable and essentially bounded input $u:\mathbb{R}_+ \to \mathbb{R}^m$ with $\|u\|_\infty \le \sigma$ the solution $x(t)$ of (5.3) is bounded, it follows that there exists a compact set $Q \subset \operatorname{int}\left(\mathbb{R}^n_+\right)$ such that for every measurable and essentially bounded input $u:\mathbb{R}_+ \to \mathbb{R}^m$ with $\|u\|_\infty \le \sigma$ the solution $y(t)$ of (5.2) is in $Q \subset \operatorname{int}\left(\mathbb{R}^n_+\right)$. Using (5.2) we get for $i=1,\dots,n$:

$$\frac{1}{T}\ln\left(\frac{y_i(T)}{y_i(0)}\right) = k_i + \sum_{j=1}^{n} B_{i,j}\left(\frac{1}{T}\int_0^T y_j(t)dt\right) + \sum_{l=1}^{m} K_{i,l}\left(\frac{1}{T}\int_0^T u_l(t)dt\right) \tag{6.38}$$

where $k=\left(k_1,\dots,k_n\right)^T$, $B=\left\{B_{i,j} : i,j=1,\dots,n\right\}$ and $K=\left\{K_{i,l} : i=1,\dots,n, l=1,\dots,m\right\}$. If $u \in C^0\left(\mathbb{R}_+;\mathbb{R}^m\right)$ is almost periodic then $\lim_{T\to+\infty}\left(\frac{1}{T}\int_0^T u(t)dt\right) \in \mathbb{R}^m$ (see [10]). Exploiting (6.38) and letting $T \to +\infty$, we get (5.5) from the fact that $k + By^* = 0$ and the fact that $y(T) \in Q \subset \operatorname{int}\left(\mathbb{R}^n_+\right)$ for all $T \ge 0$ (which implies that $\lim_{T\to+\infty}\left(\frac{1}{T}\ln\left(\frac{y_i(T)}{y_i(0)}\right)\right)=0$). The proof is complete. $\lhd$

## 7. Extensions

All the results above had to do with the existence and attractivity of almost periodic solutions. However, there is a larger class of functions that exhibits recurrence and contains all almost periodic functions: the uniformly recurrent functions (see [1, 11]). The class of uniformly recurrent functions lacks the global structural regularity of almost periodic functions and allows the study of functions with non-symmetric repetitions that can emerge from complex chaotic systems.

**Definition 4:** *A continuous function* $\varphi:\mathbb{R}\to\mathbb{R}^n$ *is called uniformly recurrent if there exists a sequence* $\{\tau_k \geq 0 : k=1,2,...\}$ *with* $\lim_{k\to+\infty}(\tau_k)=+\infty$ *such that* $\lim_{k\to+\infty}\left(\left\|\delta_{\tau_k}\varphi-\varphi\right\|_\infty\right)=0$.

**Definition 5:** *A continuous map* $f:\mathbb{R}\times\mathbb{R}^n\to\mathbb{R}^m$ *is called uniformly recurrent in* $t$ *uniformly in* $x\in\Lambda$, *where* $\Lambda\subseteq\mathbb{R}^n$ *is a non-empty set, if there exists a sequence* $\{\tau_k\geq 0 : k=1,2,...\}$ *with* $\lim_{k\to+\infty}(\tau_k)=+\infty$ *such that* $\lim_{k\to+\infty}\left(\sup_{t\in\mathbb{R},x\in\Lambda}\left(\left|f(t+\tau_k,x)-f(t,x)\right|\right)\right)=0$.

The proofs of all the results that were given above are valid with almost no change if we replace "almost periodic" by "uniformly recurrent". Thus, we obtain the following results.

**Theorem 1'-Existence of an attracting uniformly recurrent solution:** *Consider system (2.1) and assume that* $f\in C^1\left(\mathbb{R}\times\mathbb{R}^n;\mathbb{R}^n\right)$ *satisfies assumption (A). Let* $C\subseteq\mathbb{R}^n$ *be a convex, compact set for which* $f\in C^1\left(\mathbb{R}\times\mathbb{R}^n;\mathbb{R}^n\right)$ *is uniformly recurrent in* $t$ *uniformly in* $x\in C$. *Suppose that there exists* $t_0\in\mathbb{R}$, $x_0\in\mathbb{R}^n$ *and* $\tau\in\mathbb{R}$ *with* $\varphi(t,t_0,x_0)\in C$ *for all* $t\geq\tau$. *Let* $p\in C^0(\mathbb{R})$ *be a function that satisfies* $\int_0^{+\infty}p(s)ds=+\infty$ *and* $\mu\left(\frac{\partial f}{\partial z}(t,z)\right)\leq -p(t)$ *for all* $t\in\mathbb{R}$, $z\in C$, *where* $\mu(\cdot)$ *denotes the matrix measure associated with a norm* $|\cdot|$ *of* $\mathbb{R}^n$. *Suppose that there exist constants* $K,c>0$ *such that for every* $t\in\mathbb{R}$ *property (I) of Theorem 1 holds. Then there exists a uniformly recurrent solution* $x^*\in C^1(\mathbb{R};C)$ *of (2.1). Moreover,* $\lim_{t\to+\infty}\left(\left|\varphi(t,t_0,x_0)-x^*(t)\right|\right)=0$ *for all* $t_0\in\mathbb{R}$, $x_0\in S(t_0)$, *where* $S(t_0)\subseteq\mathbb{R}^n$ *is the set of all* $x_0\in\mathbb{R}^n$ *for which there exists* $\tau\in\mathbb{R}$ *with* $\varphi(t,t_0,x_0)\in C$ *for all* $t\geq\tau$.

**Lemma 1'-Conditions for the existence of a bounded entire solution:** *Let* $C\subseteq\mathbb{R}^n$ *be a compact set for which there exists* $t_0\in\mathbb{R}$, $x_0\in\mathbb{R}^n$ *and* $\tau\in\mathbb{R}$ *with* $\varphi(t,t_0,x_0)\in C$ *for all* $t\geq\tau$. *Suppose that* $f\in C^1\left(\mathbb{R}\times\mathbb{R}^n;\mathbb{R}^n\right)$ *is uniformly recurrent in* $t$ *uniformly in* $x\in C$ *and satisfies assumption (A). Then there exists an entire solution* $x^*\in C^1(\mathbb{R};C)$ *of (2.1).*

**Theorem 2'-Conditions for a bounded entire solution to be a uniformly recurrent solution:** *Consider system (2.1) and let* $C\subseteq\mathbb{R}^n$ *be a convex, compact set for which* $f\in C^1\left(\mathbb{R}\times\mathbb{R}^n;\mathbb{R}^n\right)$ *is uniformly recurrent in* $t$ *uniformly in* $x\in C$. *Let* $p,\beta\in C^0(\mathbb{R})$ *be functions that satisfy (2.2) for all* $t\in\mathbb{R}$, $z\in C$, *where* $\mu(\cdot)$ *denotes the matrix measure associated with a norm* $|\cdot|$ *of* $\mathbb{R}^n$. *Suppose that there exist constants* $K,c>0$ *such that for every* $t\in\mathbb{R}$ *either property (I) of Theorem 1 holds or property (II) of Theorem 2 holds. Furthermore, suppose that assumption (A) holds. Let* $x^*:\mathbb{R}\to C$ *be an entire solution of (2.1). Then* $x^*$ *is uniformly recurrent.*

The same observation is valid for the entrainment problem. We next give the definition of entrainment to uniformly recurrent inputs.

**Definition 6:** *Let $D \subseteq \mathbb{R}^n$ be an open set. Consider the system (4.1) with $f \in C^1\left(D \times \mathbb{R}^m; \mathbb{R}^n\right)$. We say that:*

**(a)** *system (4.1) globally entrains to uniformly recurrent inputs if the following property holds: given an input $u \in C^0\left(\mathbb{R}; \mathbb{R}^m\right)$ which is uniformly recurrent, all solutions of (4.1) converge to a unique uniformly recurrent solution.*

**(b)** *system (4.1) globally entrains to uniformly recurrent inputs of magnitude $r > 0$ if the following property holds: given an input $u \in C^0\left(\mathbb{R}; \mathbb{R}^m\right)$ with $\|u\|_\infty \le r$ which is uniformly recurrent, all solutions of (4.1) converge to a unique uniformly recurrent solution.*

**(c)** *system (4.1) globally entrains to small uniformly recurrent inputs if there exists $r > 0$ such that system (4.1) globally entrains to uniformly recurrent inputs of magnitude $r > 0$.*

The proofs of all the results for the entrainment problem that were given above are valid with almost no change if we replace “almost periodic” by “uniformly recurrent”. Thus, we obtain the following results.

**Corollary 1’:** *Consider system (4.1) with $f \in C^1\left(\mathbb{R}^n \times \mathbb{R}^m; \mathbb{R}^n\right)$ and $D = \mathbb{R}^n$. Let $|\cdot|$ denote the standard Euclidean norm of $\mathbb{R}^n$ or $\mathbb{R}^m$. Suppose that there exist constants $R, r > 0$ such that the properties (P1) and (P2) of Corollary 1 hold. Then system (4.1) globally entrains to uniformly recurrent inputs of magnitude $r$.*

**Corollary 2’:** *Consider system (4.1) with $f \in C^1\left(\mathbb{R}^n \times \mathbb{R}^m; \mathbb{R}^n\right)$ and $D = \mathbb{R}^n$. Suppose that $f\left(0,0\right) = 0$ and that system (4.1) satisfies the Input-to-State Stability (ISS) property with respect to small inputs. Moreover, suppose that $0 \in \mathbb{R}^n$ is locally exponentially stable (LES) for the input-free system $\dot{x} = f\left(x,0\right)$. Then system (4.1) globally entrains to small uniformly recurrent inputs.*

**Theorem 3’:** *Consider the Lotka-Volterra system (5.2) and suppose that there exists an interior equilibrium point $y^* \in \mathrm{int}\left(\mathbb{R}^n_+\right)$ with $k + By^* = 0$. Moreover, suppose that the matrix $A = Bdiag\left(y^*\right)$ is a Volterra-Lyapunov stable matrix, i.e., there exists $p \in \mathrm{int}\left(\mathbb{R}^n_+\right)$ such that the matrix $diag(p)A + A^T diag(p)$ is negative definite. Then system (5.2) globally entrains to small uniformly recurrent inputs. Finally, equation (5.5) holds for all uniformly recurrent inputs of sufficiently small magnitude for which the limit $\lim\limits_{T \to +\infty}\left(\frac{1}{T}\int\limits_0^T u(t)dt\right)$ exists in $\mathbb{R}^m$.*

## 8. Conclusions

The paper provided sufficient conditions for the existence and attractivity of almost periodic solutions and uniformly recurrent solutions for nonlinear time-varying systems. Moreover, the paper provided sufficient conditions for global entrainment to small almost periodic (and uniformly recurrent) inputs in nonlinear time-invariant systems. The obtained results leave room for improvement: the relaxation of the given sufficient conditions is the topic of future research.

## References

[1] E. Ait Dads, B. Es-sebbar and L. Lhachimi, "Some Constructions and Properties of Almost Automorphic Functions and Size-Structured Population Models", *Journal of Mathematical Analysis and Applications*, 562, 2026, 130731.

[2] L. Amerio and G. Prouse, *Almost Periodic Functions and Functional Equations*, Springer Verlag, 1971.

[3] F. Dörfler and F. Bullo, "Synchronization and Transient Stability in Power Networks and Nonuniform Kuramoto Oscillators", *SIAM Journal on Control and Optimization*, 50, 2012, 1616-1642.

[4] A. Duvall and E. D. Sontag, "Global Exponential Stability (and Contraction of an Unforced System) Does not Imply Entrainment to Periodic Inputs", arXiv:2310.03241 [math.OC].

[5] A. M. Fink, *Almost Periodic Differential Equations*, Lecture Notes in Mathematics, 377, Springer-Verlag, 1974.

[6] A. Genda, A. Fidlin, S. Lenci and O. Gendelman, "Effects of Seasonal Hunting and Food Quantity Variations on the Lotka-Volterra Model", *Nonlinear Dynamics*, 113, 2025, 12351–12369.

[7] M. E. Gilmore, C. Guiver and H. Logemann, "Incremental Input-to-State Stability for Lur'e Systems and Asymptotic Behaviour in the Presence of Stepanov Almost Periodic Forcing", *Journal of Differential Equations*, 300, 2021, 692–733.

[8] M. E. Gilmore, C. Guiver and H. Logemann, "Infinite-Dimensional Lur'e Systems with Almost Periodic Forcing", *Mathematics of Control, Signals, and Systems*, 32, 2020, 327–360.

[9] C. Guiver, "Passivity Theorems for Input-to-State Stability of Forced Lur'e Inclusions and Equations, and Consequent Entrainment-Type Properties", *ESAIM: Control, Optimisation and Calculus of Variations*, 31, 2025, 22.

[10] J. K. Hale, *Oscillations in Nonlinear Systems*, Dover, 1992.

[11] A. Haraux and P. Souplet, "An Example of Uniformly Recurrent Function Which is not Almost Periodic", *Journal of Fourier Analysis and Applications*, 10, 2004, 217-220.

[12] J. Hofbauer and K. Sigmund, *Evolutionary Games and Population Dynamics*, Cambridge University Press, 1998.

[13] I. Karafyllis and Z.-P. Jiang, *Stability and Stabilization of Nonlinear Systems*, Springer-Verlag, 2011.

[14] I. Karafyllis, "Results for Global Attractivity of Interior Equilibrium Points for Lotka-Volterra Systems", to appear in *Dynamical Systems* (see also arXiv:2512.11384 [math.DS]).

[15] I. Karafyllis and M. Krstic, "On the Existence of Periodic Solutions with Applications to Extremum-Seeking", submitted to *Mathematics of Control, Signals, and Systems* (see also arXiv:2602.13615 [math.OC]).

[16] H. K. Khalil, *Nonlinear Systems*, 2nd Edition, Prentice Hall, 1996.

[17] M. Krstic. "Asymmetry Demystified: Strict CLFs and Feedbacks for Predator–Prey Interconnections", arXiv:2602.21594 [eess.SY].

[18] W. Lohmiller, J.-J. E. Slotine, "On Contraction Analysis for Non-linear Systems", *Automatica*, 34, 1998, 683–696.
[19] M. Margaliot, E. D. Sontag and T. Tuller, "Entrainment to Periodic Initiation and Transition Rates in a Computational Model for Gene Translation", *PLoS ONE*, 9, 2014, e96039.
[20] M. Margaliot, E. D. Sontag and T. Tuller, "Contraction after small transients", *Automatica*, 67, 2016, 178-184.
[21] A. Mironchenko, *Input-to-State Stability. Theory and Applications*, Springer, 2023.
[22] A. Pavlov, A. Pogromsky, N. van de Wouw, H. Nijmeijer, "Convergent Dynamics, a Tribute to Boris Pavlovich Demidovich", *Systems & Control Letters*, 52, 2004, 257 − 261.
[23] G. Russo, M. di Bernardo and E. D. Sontag, "Global Entrainment of Transcriptional Systems to Periodic Inputs", *PLoS Computational Biology*, 6, 2010, e1000739.
[24] M. Vidyasagar, *Nonlinear Systems Analysis*, 2nd Edition, Prentice-Hall, 1993.

## Appendix

**Proof of (3.13):** Define the functions

$$h(z)=\exp(z)-z-1,\ \text{for}\ z\in\mathbb{R} \tag{A.1}$$

$$V(x)=h(x_1)+h(x_2),\ \text{for}\ x=(x_1,x_2)\in\mathbb{R}^2 \tag{A.2}$$

Notice that the functions $h(z)$ and $V(x)$ are both positive definite and radially unbounded. The derivative of $V(x)$ along the solutions of (3.11) satisfies:

$$\begin{aligned}
&\dot{V}(x(t))=\left(\exp(x_1(t))-1\right)\left(1-\exp(2x_1(t))\right)+a\exp(x_2(t))\left(\exp(x_1(t))-1\right)\cos(\varphi_0+\omega t)\\
&+\left(\exp(x_2(t))-1\right)\left(1-\exp(2x_2(t))\right)+b\exp(x_1(t))\left(\exp(x_2(t))-1\right)\cos(\theta_0+t)\\
&=-\left(\exp(x_1(t))+1\right)\left(\exp(x_1(t))-1\right)^2-\left(\exp(x_2(t))+1\right)\left(\exp(x_2(t))-1\right)^2\\
&+\left(\exp(x_1(t))-1\right)\left(\exp(x_2(t))-1\right)\left(a\cos(\varphi_0+\omega t)+b\cos(\theta_0+t)\right)\\
&+a\left(\exp(x_1(t))-1\right)\cos(\varphi_0+\omega t)+b\left(\exp(x_2(t))-1\right)\cos(\theta_0+t)\\
&\leq-\left(\exp(x_1(t))-1\right)^2-\left(\exp(x_2(t))-1\right)^2+\left(|a|+|b|\right)\left|\exp(x_1(t))-1\right|\left|\exp(x_2(t))-1\right|\\
&+|a|\left|\exp(x_1(t))-1\right|+|b|\left|\exp(x_2(t))-1\right|
\end{aligned} \tag{A.3}$$

Using the inequalities

$$\begin{aligned}
&\left|\exp(x_1(t))-1\right|\left|\exp(x_2(t))-1\right|\leq\frac{1}{2}\left|\exp(x_1(t))-1\right|^2+\frac{1}{2}\left|\exp(x_2(t))-1\right|^2\\
&|a|\left|\exp(x_1(t))-1\right|\leq\frac{1}{2}\left|\exp(x_1(t))-1\right|^2+\frac{a^2}{2}\\
&|b|\left|\exp(x_2(t))-1\right|\leq\frac{1}{2}\left|\exp(x_2(t))-1\right|^2+\frac{b^2}{2}
\end{aligned}$$

and since $|a|+|b|<1$ (a consequence of (3.7)), we obtain from (A.3):

$$\begin{aligned}\dot{V}\left(x(t)\right) &\leq -\frac{1-|a|-|b|}{2}\left(\left(\exp\left(x_1(t)\right)-1\right)^2+\left(\exp\left(x_2(t)\right)-1\right)^2\right)+\frac{a^2+b^2}{2} \\ &\leq -\frac{1-|a|-|b|}{2}\left(\max_{i=1,2}\left(\left|\exp\left(x_i(t)\right)-1\right|\right)\right)^2+\frac{a^2+b^2}{2}\end{aligned} \tag{A.4}$$

The following inequality holds for all $z \in \mathbb{R}$

$$\left|\exp(z)-1\right| \geq \min\left(\exp(z)-z-1, \frac{1}{2}\right) \tag{A.5}$$

and since $|a|+|b|<1$, we obtain from (A.1) and (A.4):

$$\begin{aligned}\dot{V}\left(x(t)\right) &\leq -\frac{1-|a|-|b|}{2}\left(\max_{i=1,2}\left(\min\left(h\left(x_i(t)\right), \frac{1}{2}\right)\right)\right)^2+\frac{a^2+b^2}{2} \\ &= -\frac{1-|a|-|b|}{2}\left(\min\left(\max_{i=1,2}\left(h\left(x_i(t)\right)\right), \frac{1}{2}\right)\right)^2+\frac{a^2+b^2}{2}\end{aligned} \tag{A.6}$$

Using definition (A.2) we get the following inequality

$$V(x) \leq 2\max_{i=1,2}\left(h\left(x_i\right)\right) \tag{A.7}$$

which combined with (A.6) and the fact that $|a|+|b|<1$, gives:

$$\dot{V}\left(x(t)\right) \leq -\frac{1-|a|-|b|}{8}\left(\min\left(V\left(x(t)\right),1\right)\right)^2+\frac{a^2+b^2}{2} \tag{A.8}$$

The differential inequality (A.8) and (3.7) (which implies that $-\frac{1-|a|-|b|}{8}+\frac{a^2+b^2}{2}<0$) guarantee that for every $\lambda \in (0,1)$ which is sufficiently close to 1 and so that $-\lambda\frac{1-|a|-|b|}{8}+\frac{a^2+b^2}{2}\leq 0$, the following implication holds:

$$V\left(x(t)\right) \geq 2\sqrt{\frac{a^2+b^2}{\lambda\left(1-|a|-|b|\right)}} \Rightarrow \dot{V}\left(x(t)\right) \leq -(1-\lambda)\frac{1-|a|-|b|}{8}\left(\min\left(V\left(x(t)\right),1\right)\right)^2 \tag{A.9}$$

Implication (A.9) in conjunction with Lemma 2.14 on page 82 in [13] and the fact that $V(x)$ is positive definite and radially unbounded guarantees that every solution of (3.11) is defined for all $t \geq t_0$ and satisfies the following asymptotic estimate:

$$\limsup_{t\to+\infty}\left(V\left(x(t)\right)\right)\leq 2\sqrt{\frac{a^2+b^2}{\lambda\left(1-|a|-|b|\right)}} \tag{A.10}$$

By letting $\lambda \to 1^-$ in estimate (A.10), we get the estimate:

$$\limsup_{t\to+\infty}\left(V\left(x(t)\right)\right)\leq 2\sqrt{\frac{a^2+b^2}{1-|a|-|b|}} \tag{A.11}$$

It is clear from (A.2) and (A.11) that the following estimates hold for $i=1,2$:

$$\limsup_{t\to+\infty}\left(h\left(x_i(t)\right)\right)\leq 2\sqrt{\frac{a^2+b^2}{1-|a|-|b|}} \tag{A.12}$$

We next show that the following implication is valid for all $z\in\mathbb{R}$:

$$|z|\leq \kappa\left(\exp\left(z\right)-z-1\right) \tag{A.13}$$

where $\kappa\left(s\right)=s+\sqrt{2s}$ for $s\geq 0$.

Set $s=\exp\left(z\right)-z-1$. If $s=0$ then $z=0$ and (A.13) holds. So, next we assume that $s=\exp\left(z\right)-z-1>0$. Due to the monotonicity of $h(z)=\exp(z)-z-1$ (increasing for $z\geq 0$ and decreasing for $z\leq 0$) the equation $s=\exp\left(z\right)-z-1$ has exactly two solutions: one solution $z_1>0$ and one solution $z_2<0$. Consequently, $|z|\leq \max\left(z_1,-z_2\right)$. For $z_1>0$ we get $s=\exp\left(z_1\right)-z_1-1\geq \frac{1}{2}z_1^2$ which implies $z_1\leq\sqrt{2s}$. For $z_2<0$ we get $h\left(-s-\sqrt{2s}\right)\geq h\left(z_2\right)=s$ (the inequality $e^{-s-\sqrt{2s}}+s+\sqrt{2s}-1\geq s$ is implied by the fact that the function $f(s)=e^{-s-\sqrt{2s}}+\sqrt{2s}$ for $s\geq 0$ is non-decreasing due to the fact that $f'(s)=\frac{e^{-s-\sqrt{2s}}}{\sqrt{2s}}\left(e^{s+\sqrt{2s}}-1-\sqrt{2s}\right)\geq 0$ for $s>0$). Consequently (since $h(z)=\exp(z)-z-1$ is decreasing for $z\leq 0$), we get $-s-\sqrt{2s}\leq z_2$ which gives $-z_2\leq s+\sqrt{2s}$. Combining the facts that $s=\exp\left(z\right)-z-1>0$, $|z|\leq\max\left(z_1,-z_2\right)$, $z_1\leq\sqrt{2s}$ and $-z_2\leq s+\sqrt{2s}$, we obtain (A.13).

Exploiting (A.1), (A.12) and (A.13) we obtain (3.13). The proof is complete. $\lhd$